\documentclass[letterarticle,12pt,peerreviewca,draftcls]{IEEEtran}
\usepackage{csm16}
\usepackage[margin=1in]{geometry}
\usepackage[nolists,nomarkers,tablesfirst,disable]{endfloat} 
\usepackage{amsmath} 
\usepackage[nomarkers]{endfloat}
\usepackage{url}
\usepackage{graphicx,xcolor}
\usepackage{verbatim}
\usepackage[left,pagewise]{lineno}
\usepackage[english]{babel}
\usepackage{blindtext}
\usepackage{amssymb,amsfonts,amscd,mathrsfs}
\usepackage{upgreek,dsfont}
\usepackage{enumerate}

\makeatletter
\newcommand{\verbatimfont}[1]{\def\verbatim@font{#1}}%
\makeatother
\verbatimfont{\ttfamily\small}

\title{Dissipativity and Integral Quadratic Constraints\\
\Large Tailored computational robustness tests for complex interconnections 

\footnote{\small
\textcopyright~
20XX IEEE. Personal use of this material is permitted. Permission from IEEE must be obtained
for all other uses, in any current or future media, including reprinting/republishing this material for
advertising or promotional purposes, creating new collective works, for resale or redistribution to servers or
lists, or reuse of any copyrighted component of this work in other works. {\tt https://doi.org/10.1109/MCS.2022.3157117}}
}
\author{Carsten W. Scherer}

\newif\ifPDF \ifx\pdfoutput\undefined\PDFfalse \else\ifnum\pdfoutput > 0\PDFtrue \else\PDFfalse \fi \fi
\ifPDF
\usepackage[pdftex, plainpages = false, colorlinks=true, linkcolor=black, citecolor = green!50!blue, urlcolor = blue, filecolor=black, pagebackref=false, hypertexnames=false,  pdfpagelabels ]{hyperref}
\fi

\newcommand{\mstrut}[1]{\rule{0pt}{#1}}
\let\oldhline=\hline 
\renewcommand{\hline}{\oldhline\mstrut{2.1ex}}

\renewcommand{\S}{\Bbb S}

\newcommand{\gai}{\rm br}
\newcommand{\pas}{\rm pr}

\newcommand{\nz}{k}
\newcommand{\nw}{m}
\newcommand{\ny}{k}
\newcommand{\nuu}{m}
\newcommand{\nx}{{n}}

\newcommand{\ve}{\varepsilon}

\newcommand{\nl}{\vp}
\newcommand{\sat}{\operatorname{sat}}
\newcommand{\sign}{\operatorname{sign}}

\newcommand{\sy}{G}
\newcommand{\mi}{P}

\renewcommand{\L}{{\mathscr{L}}}

\newcommand{\de}{\delta}

\newcommand{\col}{\text{col}}
\newcommand{\pl}{\bul}
\newcommand{\Le}{{\cal L}}

\newcommand{\Se}{{\cal S}}
\newcommand{\Ri}{\Rset\cup\{\infty\}}

\newcommand{\Xs}{X}

\newcommand{\st}{\ |\ }

\newcommand{\eps}{\ve}

\newcommand{\skipthis}[1]{}
\newcommand{\bm}[1]{ \mbox{\boldmath$ #1 $} }

\renewcommand{\t}[1]{ \tilde{#1} }

\renewcommand{\c}[1]{{\cal #1}}

\newcommand{\epro}{\hfill\mbox{\rule{2mm}{2mm}}}

\newcommand{\R}{ {\mathbb{R}} }
\newcommand{\Rset}{ {\mathbb{R}} }

\newcommand{\N}{ {\mathbb{N}} }

\newcommand{\diag}{ \operatornamewithlimits{diag} }

\renewcommand{\r}[1]{(\ref{#1})}

\newtheorem{theo}{Theorem}

\newtheorem{hypo}[theo]{Hypothesis}
\newtheorem{lemm}[theo]{Lemma}
\newtheorem{coro}[theo]{Corollary}
\newtheorem{defi}[theo]{Definition}
\newtheorem{rema}[theo]{Remark}
\newtheorem{exam}[theo]{Example}
\newtheorem{assumption}[theo]{Assumption}

\newcommand{\theorem}[1]{\begin{theo}#1\end{theo}}
\newcommand{\corollary}[1]{\begin{coro}#1\end{coro}}
\newcommand{\lemma}[1]{\begin{lemm}#1\end{lemm}}
\newcommand{\definition}[1]{\begin{defi}#1\end{defi}}
\newcommand{\example}[1]{\begin{exam}#1\end{exam}}

\newcommand{\proo}[1]{{\em Proof:} #1\epro}

\newcommand{\mun}[1]{\begin{multline*}#1\end{multline*}}
\newcommand{\mul}[1]{\begin{multline}#1\end{multline}}
\newcommand{\equ}[1]{\begin{equation}#1\end{equation}}
\newcommand{\eqn}[1]{$$#1$$}
\newcommand{\eql}[2]{\begin{equation}\label{#1}#2\end{equation}}

\newcommand{\gan}[1]{\begin{gather*}#1\end{gather*}}
\newcommand{\ear}[1]{\begin{eqnarray}#1\end{eqnarr{{{}}}ay}}

\newcommand{\arr}[2]{\begin{array}{#1}#2\end{array}}
\newcommand{\mas}[2]{\left[\begin{array}{#1}#2\end{array}\right]}
\newcommand{\mat}[2]{\left(\begin{array}{#1}#2\end{array}\right)}

\newcommand{\enu}[1]{\begin{enumerate}\itemsep-1ex#1\end{enumerate}}
\newcommand{\cen}[1]{\begin{center}#1\end{center}}

\newcommand{\la}{\lambda}
\newcommand{\ga}{\gamma}
\newcommand{\om}{\omega}
\newcommand{\io}{i\omega}
\newcommand{\te}[1]{\text{\ \ #1\ \ }}
\newcommand{\be}{\beta}

\newcommand{\al}{\alpha}

\newcommand{\cl}{\prec}
\newcommand{\cg}{\succ}
\newcommand{\cge}{\succeq}
\newcommand{\cle}{\preceq}

\newcommand{\hl}{\\\hline}

\newcommand{\remark}[1]{\begin{rema}#1\end{rema}}

\usepackage{tikz,pgfplots,datatool}
\usetikzlibrary{positioning,calc,arrows,shapes,shadows}
\tikzset{
auto,
sysw/.style 2 args={
rectangle,
draw,
rounded corners,
draw=black,
fill=white,
minimum height=#2,
minimum width=#1,
inner sep=\dn},
sys/.style 2 args={
rectangle,
draw,
rounded corners,
drop shadow,
fill=white,
minimum height=#2,
minimum width=#1,
inner sep=\dn},
sum/.style={circle,draw,draw=black,inner sep=0mm,minimum size=2mm},
jun/.style={circle,draw,draw=black,inner sep=0mm,minimum size=0.5mm,fill},
>={latex},
every path/.style={rounded corners},
}

\newcommand{\tio}[4]{\coordinate (#1) at ($(#2.south #3)!#4!(#2.north #3)$)}

\newcommand{\tsh}[2]{([xshift=#2] #1)}
\def\dn{1ex}
\def\dl{3*\dn}
\def\dh{6*\dn}
\tikzstyle{sy0}=[sys={0*\dn}{0*\dn}]
\tikzstyle{sy1}=[sys={12*\dn}{8*\dn}]
\tikzstyle{sy2}=[sys={8*\dn}{6*\dn}]
\tikzstyle{sy3}=[sys={5*\dn}{5*\dn}]
\tikzstyle{sye0}=[syse={0*\dn}{0*\dn}]
\tikzstyle{sye1}=[syse={12*\dn}{8*\dn}]
\tikzstyle{sye2}=[syse={8*\dn}{6*\dn}]
\tikzstyle{sye3}=[syse={5*\dn}{5*\dn}]
\tikzstyle{sy4}=[sys={12*\dn}{6*\dn}]
\tikzstyle{sye4}=[syse={12*\dn}{6*\dn}]

\newcommand{\C}{\mathbb C}
\newcommand{\eig}{\text{eig}}
\newcommand{\bul}{\bullet}
\renewcommand{\Re}{\text{Re}}
\renewcommand{\Im}{\text{Im}}

\newcommand{\Del}{\Delta}
\newcommand{\Des}{\bm{\Delta}}
\newcommand{\Denl}{ { \bm{\Delta_{\rm nl}} } }
\newcommand{\Dednl}{ { \bm{\Delta_{\rm dnl}} } }
\newcommand{\Dedy}{ { \bm{\Delta_{\rm dy}} } }

\usepackage[compress]{cite}
\usepackage{array}
\usepackage{boldline}

\newcolumntype{?}{!{\vrule width 1pt}}

\newcommand{\vp}{\varphi}
\newcommand{\fbs}[1]{
\begin{tikzpicture}[xscale=1,yscale=1]
\tikzstyle{sy4}=[sys={5*\dn}{4*\dn}]
\def\ds{7*\dn}
\node[sy4] (g) at (0,0)  {$\sy$};
\node[sy4] (d) at ( $(g)+(0,\dh)$ ) {#1};
\draw[<-] (g)-- ($(g) + (\dh, 0)$)  |- node[swap]{} ($(g) + (\dh, \dl)$) |- (d);
\draw[->] (g)-- ($(g) + (-\dh, 0)$)  |- node[]{} ($(g) + (-\dh, \dl)$) |- (d);
\end{tikzpicture}
}

\newcommand{\fstaba}[1]{
\begin{tikzpicture}[xscale=1,yscale=1]
\tikzstyle{sy5}=[sys={5*\dn}{4*\dn}]
\def\ds{7*\dn}
\node[sy5] (g) at (0,0)  {$\sy$};
\node[sy5] (d) at ( $(g)+(0,\dh)$ ) {#1};
\node[sum]  (sum) at ($(g) + (-\ds, 0)$) {\scalebox{0.7}{$+$}};
\coordinate (jun) at ($(g) + (\ds, 0)$) {};

\draw[->] (g) -- 
 (sum);
\draw[->] (sum)   |-  node[pos=.75]{$z$}  (d);
\draw[-] (d) -| ($(g)+(\dh,0)$)  -- (g);
\draw[->] ($(sum)+(-4*\dn,0)$) -- node[]{$d$} (sum);
\end{tikzpicture}
}

\newcommand{\fstabb}[1]{
\begin{tikzpicture}[xscale=1,yscale=1]
\tikzstyle{sy5}=[sys={5*\dn}{4*\dn}]
\def\ds{7*\dn}
\node[sy5] (g) at (0,0)  {$\sy$};
\node[sy5] (d) at ( $(g)+(0,\dh)$ ) {#1};
\node[sum]  (sum) at ($(g) + (-\ds, 0)$) {\scalebox{0.7}{$+$}};
\node[jun]  (jun) at ($(g) + (\ds, 0)$) {};

\draw[->] (g) -- 
 (sum);
\draw[->] (sum)   |-  node[pos=.75]{$z$}  (d);
\draw[-] (d)  -| node[pos=.25]{} (jun);
\draw[->] (jun)   --  (g);

\draw[->] ($(sum)+(-4*\dn,0)$) -- node[]{$d$} (sum);
\draw[->] (jun)-- node{$w$} ($(jun)+(4*\dn,0)$);
\end{tikzpicture}
}

\newcommand{\fstabc}[1]{
\begin{tikzpicture}[xscale=1,yscale=1]
\tikzstyle{sy5}=[sys={5*\dn}{4*\dn}]
\def\ds{7*\dn}
\node[sy5] (g) at (0,0)  {$\sy$};
\node[sy5] (d) at ( $(g)+(0,\dh)$ ) {#1};
\node[sum]  (sum) at ($(g) + (-\ds, 0)$) {\scalebox{0.7}{$+$}};
\node[sum]  (sum2) at ($(d) + (\ds, 0)$) {\scalebox{0.7}{$+$}};

\draw[->] (g) -- 
 (sum);
\draw[->] (sum)   |-  node[pos=.75]{$z$}  (d);
\draw[-] (d)  -- (sum2);
\draw[-] (sum2) |- node[pos=.75,swap]{$u$} (g);
\draw[->] ($(sum)+(-4*\dn,0)$) -- node[]{$d_1$} (sum);
\draw[->] ($(sum2)+(4*\dn,0)$) -- node[swap]{$d_2$} (sum2);
\end{tikzpicture}
}

\newcommand{\fbp}[1]{
\begin{tikzpicture}
\def\dnl{2*\dn}
\node[sy2] (g) at (0,.3) {$\mat{cc}{\sy&\sy_{d}\\\sy_e&\sy_{ed}}$};
\tio{i1}{g}{west}{2/3};
\tio{o1}{g}{east}{2/3};
\tio{i2}{g}{west}{1/3};
\tio{o2}{g}{east}{1/3};
\node[sy3, above = \dnl of g] (d) {#1};
\draw[->] (i1) -- \tsh{i1}{-1*\dl}  |- node[pos=.25] {$z$} (d);
\draw[<-] (o1) -- \tsh{o1}{1*\dl}|- node[pos=.25,swap] {$w$} (d);
\draw[<-] (o2) -- node[swap]{$d$} \tsh{o2}{1*\dh};
\draw[->] (i2) -- node[pos=.5]{$e$} \tsh{i2}{-1*\dh};

\tikzstyle{sy4}=[sys={5*\dn}{4*\dn}]
\def\ds{7*\dn}
\node[sy4] (g) at (-12,0)  {$\sy$};
\node[sy4] (d) at ( $(g)+(0,\dh)$ ) {#1};
\draw[<-] (g)-- ($(g) + (\dh, 0)$)  |- node[swap]{$w$} ($(g) + (\dh, \dl)$) |- (d);
\draw[->] (g)-- ($(g) + (-\dh, 0)$)  |- node[]{$z$} ($(g) + (-\dh, \dl)$) |- (d);

\node[sy4] (g) at (-7,0)  {$\sy$};
\node[sy4] (d) at ( $(g)+(0,\dh)$ ) {#1};
\draw[<-] (g)-- ($(g) + (1.1*\dh, 0)$)  |- node[swap]{$w$} ($(g) + (1.1*\dh, \dl)$) |- (d);
\draw[->] (g)-- ($(g) + (-1.1*\dh, 0)$)  |- node[]{$z$} ($(g) + (-1.1*\dh, \dl)$) |- (d);
\tio{Di1}{d}{west}{1/2};
\tio{Do1}{d}{east}{1/2};
\node[sys={4*\dn}{5*\dn}] (f) at ($(d)+(2,.7)$)  {$\Psi$};
\tio{i1}{f}{west}{1/3};
\tio{i2}{f}{west}{2/3};
\draw[->] ([xshift=.7*\dl] Do1) |- (i1);
\draw[->] ([xshift=-.7*\dl] Di1) |- (i2);
\node[jun] (jun1) at ([xshift=-.7*\dl] Di1) {};
\node[jun] (jun2) at ([xshift=.7*\dl] Do1) {};
\draw[->] (f)--node[]{$v$} ([xshift=2*\dl] f.east) ;
\end{tikzpicture}
}

\usepackage[absolute]{textpos}
\begin{document}

\maketitle

\CSMsetup

\begin{textblock}{13}(1.15, 16)
\fbox{
\textcopyright~ 
20XX IEEE. Personal use of this material is permitted. Permission from IEEE must be obtained
for all other uses, in any current or future media, including reprinting/republishing this material for
advertising or promotional purposes, creating new collective works, for resale or redistribution to servers or
lists, or reuse of any copyrighted component of this work in other works.}
\end{textblock}

The field of systems and control is dominated by understanding complex cyberphysical dynamical systems or designing such systems to achieve a desired dynamic behavior. To address complexity, it is key 
to consider a monolithic system as an interconnection of subsystems. This permits capturing dynamical properties of systems at the interconnection level by the characteristics of the subsystems and  their interconnection.

Analyzing the stability of a feedback loop of two systems is a ubiquitous 
scenario of this kind.
It is a natural question to ask which properties of the individual systems guarantee interconnection stability.
The classical passivity and small-gain theorems are incarnations of such a strategy. Myriads of related so-called absolute stability results have been proposed in the control literature since the formalization of these concepts in the $1960$s. Jan Willems' introduction in \cite{Bas01} to the seminal article \cite{Pop62}, which emphasizes the substantial impact of Popov's work on the development of absolute stability theory (for example, in \cite{Yak62,BroWil65a,Zam66,Zam66b}). In commenting on Zames' foundational contribution
\cite{Zam66b}, Jan
notes  in \cite{Bas01}
that this work anticipates the power of multipliers to render the Popov criterion more flexible and less conservative.
Beautiful examples are the use of so-called Zames-Falb multipliers in \cite{Osh67,WilBro68,ZamFal68} (see also
\cite{Wil71book,DesVid75}). Today, absolute stability theory is considered an essential ingredient of feedback stability analysis
and has had a striking impact on the development of robust and $H_\infty$-control.

Willems developed dissipativity theory as an extension of Lyapunov theory.
This allowed many of the ideas available for closed dynamical systems (which are not influenced by the environment in which
they operate) to open dynamical systems that have inputs and outputs.
It was his explicit goal to arrive at a unifying view of
the stability properties of general feedback interconnections \cite{Wil72a,Wil72b}.

The following quote from the
the third part \cite{Wil71}  of his trilogy anticipates the impact of the development of solvers \cite{NesNem94,BoyVan04} for linear matrix inequality problems onto the field \cite{BoyGha94}:
\begin{quote}\small
``The basic importance of the LMI seems to be largely unappreciated. It would be interesting to see whether or not it can be exploited in computational algorithms, for example.''
\end{quote}

Despite the progress concerning the computation of multipliers or scalings in structured singular value theory \cite{DoyPac91} in the 1980s,  it has been lamented in \cite{MegRan97} that the use of numerical algorithms in classical absolute stability theory is limited through causality requirements on multipliers. Building on the seminal contributions of Yakubovich \cite{Yak67},
the framework of integral quadratic constraints (IQCs) was developed in \cite{MegRan97} to
analyze the stability of an interconnection of some linear time-invariant system in feedback with another causal system without any particular description (also called uncertainty). This framework does indeed offer a modular encompassing computational approach for the stability analysis of structured interconnections. Its stunningly wide impact incorporates, among many others, the analysis of adaptive learning \cite{AndYou07}
or quadratic \cite{HeaWil07} and first-order convex optimization algorithms \cite{LesRec16}.
However, the underlying technique of proof is functional analytic and, unlike with dissipativity,
does not directly allow for drawing conclusions for the state of the system in the loop, for example.

When Siep Weiland and I set out to write the first versions of our lecture notes on linear matrix inequalities in control about 25 years ago \cite{SchWei99}, it was an obvious choice to cover both dissipativity and IQC theory. One of our goals was to prove the main result of IQC theory by dissipativity arguments to fully exploit the benefits of both worlds.
This was successful for so-called hard dynamic IQCs. However, we failed to completely cover genuine soft versions thereof. 

Since then, a rather large body of work has been devoted to analyzing the links between both frameworks, and we feel unable to do justice to all these developments in the current article (with more references appearing in \cite{SchVee18}).
Still, for the genuine framework from \cite{MegRan97}, the connection between dissipativity and IQC theory was clarified
only rather recently for specialized cases in \cite{Bal02,WilTak07,VeeSch13a,Sei15,CarSei18}.
In the author's view, a complete understanding of this link can be based on the notion of
IQCs with terminal cost, as introduced in \cite{SchVee18} and further developed here.

The article is structured as follows. In the section ``\nameref{S1},'' the setup for robust stability analysis is motivated by recalling the classical Lur'e problem. This section comprises a formulation of the circle criterion, the proof of which is an ideal illustration for how dissipativity arguments are leveraged for robustness analysis. A subsequent section collects several essential results about dissipativity,
with an emphasis on linear systems and quadratic supply rates. This part exposes the classical cycle of proof, relating dissipation inequalities, frequency domain inequalities, and linear matrix inequalities (LMIs), while also covering strict versions thereof under minimal assumptions.
The ``\nameref{sb:dis}'' section is provided for the convenience of the reader.
The section ``\nameref{RS}'' covers the key technical framework for robustness analysis based on the notion of IQCs with terminal cost. Proofs are dissipation-based, while
``\nameref{sb:iqc}'' exposes in full detail why the proposed framework seamlessly encompasses
dissipativity and IQC approaches. A final section of the article shows how to expand the main robust stability and performance results if the uncertainties are structured. It is well illustrated that structure is key to both capturing complex system interconnections and the construction of tailored computational robustness tests. Novel results are provided for dynamic linear time-invariant uncertainties, which are accompanied by a numerical illustration.

In a nutshell, the current  article serves to unveil the principles underlying the systematic construction of robustness tests for structured interconnections that are amenable to computations using any
 LMI solver engine.
This reduces to paraphrasing Willems' central result about the dissipativity of a neutral interconnection of dissipative subsystems in \cite{Wil72a}, in combination with his later work on tearing, zooming, and linking
\cite{Wil07csm}, and his vision that LMIs could play a major role in computations \cite{Wil71}.

\section[Motivation: The Lur'e problem]{Motivation: The Lur'e problem}\label{S1}

A large part of systems theory has been devoted to the analysis of stability of feedback systems.
A prominent example is the Lur'e problem, which guarantees the stability of a linear system $\sy$ in feedback with a static nonlinearity $\nl$, as shown on the left in Fig.~\ref{fig1} \cite{AizGan64,JayLog11}. Such a configuration might
result from a standard saturated control loop with a linear plant $G$ and controller $K$ for $d=0$, which is depicted on the right in
Fig.~\ref{fig1}. However, it also covers many others such as deep neural networks \cite{SuyVan99}.

\subsection{Responses of feedback interconnections}

More precisely, the lower block $\sy$ in Fig.~\ref{fig1} is a linear system that admits the description
\eql{sy}{\dot x=Ax+Bu,\ \ y=Cx+Du}
with given matrices $A\in\R^{\nx\times\nx}$, $B\in\R^{\nx\times\nuu}$, $C\in\R^{\ny\times\nx}$,
and $D\in\R^{\ny\times\nuu}$. The transfer matrix of \r{sy} consists of real rational proper transfer functions as its entries and is defined by
\eql{tf}{\sy(s)=C(sI-A)^{-1}B+D,\te{which is also expressed as}\sy=\mas{c?c}{A&B\\\hlineB{2} C&D}.}
The upper block in Fig.~\ref{fig1} is a static nonlinear system defined as
\eql{vp}{
w(t)=\nl(z(t))\te{for}t\geq 0\te{with some map}\nl:\R^\nz\to\R^\nw.
}
According to the block diagram in Fig.~\ref{fig1}, these two systems are interconnected by signal sharing \cite{Wil07csm}, by equating inputs and outputs according to the relations $u=w$, and $z=y$.
This leads to the description of the feedback loop as
\eql{fb}{\dot x=Ax+Bw,\ \ z=Cx+Dw,\ \ w=\nl(z),\te{and}x(0)=x_0.}
Here, the initial condition $x_0\in\R^\nx$ is viewed as an externally specified excitation or disturbance.

All systems in this article are described on the positive real line $[0,\infty)$ as their time set.
For mathematical preciseness, we need to specify the regularity of signals as functions of time.
We assume that any free signal, such as an input or a disturbance, is taken from the universe
$$
\Se^\nuu:=\{u:[0,\infty)\to\R^\nuu\st u\text{\ is piece-wise continuous and right-continuous}\}.
$$
If $u\in\Se^\nuu$ is square integrable, its $\Le_2$-norm is denoted as $\|u\|_2:=\sqrt{\int_0^\infty u(t)^Tu(t)\,dt}$.
Moreover, state trajectories $x:[0,\infty)\to\R^\nx$ of systems described by linear or nonlinear differential equations are required to be continuous and have a right-derivative $\dot x(t)$ for each $t\geq 0$ with $\dot x\in\Se^n$.
For example \cite{LogRya14} provides an exposition of all technical details about the theory
of ordinary differential equations (ODEs) required here.

Given $u\in\Se^\nuu$ and an initial condition $x_0\in\R^\nx$, by standard ODE theory, \r{sy} admits a unique state and output response that satisfy $x(0)=x_0$ and
\r{sy} everywhere on the time set.
Similarly, the system defined by \r{vp} with a continuous map $\vp$ takes signals $z\in\Se^\nz$ into $w\in\Se^\nw$.
However, $\nl$ can be more general and admit multiple values. Then, the system defined by $\vp$ is the set of all pairs $(z,w)\in\Se^\nz\times\Se^\nw$ for which \r{vp} is valid. This follows \cite{Pop62,Zam66} in classical stability theory and is very much in the spirit of Willems' behavioral approach to systems theory \cite{Wil07csm}. In this article,   there is no need for any further formalization without harm.
The same viewpoint is taken for the feedback interconnection \r{fb}.
Any triple  $(x,w,z)\in\Se^\nx\times\Se^\nw\times\Se^\nz$ of signals that satisfies \r{fb}
is called a response of the feedback loop to the disturbance $x_0\in\R^\nx$.
Even if the initial condition $x_0$ is fixed, there might exist no, exactly one, or a set of responses of the loop.

If $D$ vanishes, \r{fb} can be equivalently expressed as
\eql{fbe}{\dot x=Ax+B\nl(Cx),\ \ z=Cx,\ \ w=\nl(Cx)\te{and}x(0)=x_0.}
Then standard ODE theory identifies properties of $\nl$ such that the initial value problem in \r{fbe} has a unique solution $x\in\Se^\nx$ that exists on the whole time set $[0,\infty)$. This is called the state response of the loop, and the remaining equations uniquely define the other signals $u=z$ and $y=w$. This well-posedness property is true if, for example,  $\nl$ is locally Lipschitz and linearly bounded \cite{LogRya14}. However, it is emphasized that no such assumptions about the existence and uniqueness of responses of differential equations or feedback loops are required
in this article.

\subsection{Definition of stability}\label{Sstab}

Many notions of stability for \r{fb} have been studied in the literature.  We follow a classical path and require the existence of some $\ga>0$ such that all responses of \r{fb} satisfy
\eql{sta}{
\int_0^T\|x(t)\|^2\,dt+\int_0^T\|w(t)\|^2\,dt\leq \ga^2\|x_0\|^2\te{for all}T\geq 0\te{and all}x_0\in\R^\nx.
}
This guarantees that $x$ and $w$ are of finite energy and that $\sqrt{\|x\|_2^2+\|w\|_2^2}\leq\ga\|x_0\|$ holds for all loop trajectories, in accordance with the classical notion of the gain of a system in \cite{Zam66}.
Moreover, from $\dot x=Ax+Bw$, infer that $\dot x$ has finite energy with
$\|\dot x\|_2\leq \ga\|(A\ B)\|\|x_0\|$, where $\|.\|$ denotes both the Euclidean norm and the spectral norm for real and complex vectors and matrices. Then, $\|x(t)\|^2=\|x_0\|^2+\int_0^t 2x(\tau)^T\dot x(\tau)\,dt$ for $t\geq 0$ shows with $K:=\sqrt{1+2\ga^2\|(A\ B)\|}$ that
$\lim_{t\to\infty}x(t)=0$ and $\sup_{t\geq 0}\|x(t)\|\leq K\|x_0\|$ holds for all loop trajectories.
This assures global asymptotic stability in the sense of Lyapunov.
Further variants of stability characterizations
in relation to the Lur'e problem are collected in \cite{JayLog11}.

\subsection{Sector bounded nonlinearities}

In the Lur'e problem, the nonlinearity $\nl$ is assumed to satisfy the  sector condition
\eql{sec}{
(Lz-\nl(z))^T(\nl(z)-Mz)\geq 0\te{for all}z\in\R^\nz
}
with two fixed matrices $M,L\in\R^{\nw\times\nz}$. The naming is motivated by considering $\nw=\nz=1$ and $M\leq L$.
Then, \r{sec}  means that the graph of $\nl$ is located in
the  conic sector defined by the graphs of the linear functions $\nl_L(z):=Lz$  and $\nl_M(z):=Mz$.
Much of classical absolute stability theory is motivated by
considering the saturation function defined as $\sat(z)=z$ for $z\in[-1,1]$ and $\sat(z)=\sign(z)$ for $|z|>1$, which satisfies a sector bound with $M=0$ and $L=1$.

With a real symmetric or complex Hermitian matrix $P$ of dimension $\ny+\nuu$
(and if $\mbox{}^*$ denotes conjugate transposition), introduce the abbreviation
$$
s_P(U,Y):=\mat{c}{Y\\U}^*\!\!P\mat{c}{Y\\U}\in\C^{p\times p}\te{for}(U,Y)\in\C^{\nuu\times p} \times\C^{\ny\times p}.
$$
Then, $(Lz-\nl(z))^T(\nl(z)-Mz)=\frac{1}{2}s_{P_{L,M}}(\vp(z),z)$ holds with the real symmetric matrix
$$
P_{L,M}:=\mat{cc}{L&-I\\-M&I}^T\!\!\mat{cc}{0&I\\I&0}\mat{cc}{L&-I\\-M&I}=
\mat{cc}{-L^TM-M^TL&L^T+M^T\\L+M&-2I}.
$$
Therefore, \r{sec} is a quadratic constraint on the graph of $\vp$ and reads with $P=P_{L,M}$ as
\eql{qccir}{
s_P(\nl(z),z)=
\mat{c}{z\\\nl(z)}^T\!\!P\mat{c}{z\\\nl(z)}\geq 0\te{for all}z\in\R^\nz.}

\subsection{The circle criterion}

Stability criteria serve to impose conditions on \r{sy} such that the feedback loop \r{fb} is stable if $\nl$ satisfies \r{qccir}. We are confronted with a question of robust or, equivalently, absolute stability against the uncertainty comprising all maps $\nl:\R^\nz\to\R^\nw$ with \r{qccir}.
For $P=P_{L,M}$, this quadratic constraint is a consequence of a sector condition. However, the following results
require no particular structure of $P\in\S^{\ny+\nuu}$, where $\S^n$ denotes the subset of symmetric matrices in $\R^{n\times n}$.

All robust stability tests involve some assumption on nominal stability, which means to assume stability for one element in the
uncertainty set. Specifically, it is assumed that \r{fb} is stable for some linear map $\nl$ described as $\vp(z)=\Del_0z$, where $\Delta_0$ satisfies
\eql{iqcnom}{
\mat{c}{I\\\Del_0}^T\!\!P\mat{c}{I\\\Del_0}\cge 0\te{and}\Del_0\in\R^{\nuu\times\ny}.
}
Note that $M\cge 0$ (or $M\cg 0$) means that the matrix $M$ is real symmetric or complex Hermitian and
positive semidefinite (or positive definite).
This leads to the following version of the circle criterion in terms of a so-called frequency-domain inequality (FDI) imposed on the transfer matrix \r{tf} of the linear system \r{sy}.

\theorem{\label{Tcir}Let $A$ in \r{sy} have no eigenvalues in $i\R$. Suppose that there exists some $\Del_0$ with \r{iqcnom} such that
all trajectories of the loop \r{fb} for $\nl(z)=\Del_0z$ satisfy $\lim_{t\to\infty}x(t)=0$.
Then \r{fb} is stable for all nonlinearities $\nl$ satisfying \r{qccir} if the following strict FDI holds:
\eql{fdicir}{
\mat{c}{\sy(i\om)\\I}^*\!\!P\mat{c}{\sy(i\om)\\I}\cl 0
\te{for all}\om\in \R\cup\{\infty\}.}}

The proof anticipates the interpretation of both \r{qccir} and \r{fdicir} as dissipation inequalities (Theorem~\ref{Tddis} and Corollary~\ref{Csdi}) and is presented
in ``\nameref{sb:pcir}.'' In particular, it shows that
\r{iqcnom} and \r{fdicir} for $\om=\infty$ imply that $I-D\Del_0$ is invertible.
For $\vp(z)=\Delta_0z$, the feedback loop \r{fb} can then also be described as
\eql{fbl}{\dot x=(A+B\Del_0(I-D\Del_0)^{-1}C)x,\ \ z=(I-D\Del_0)^{-1}Cx,\ \ w=\Del_0 z,\te{and}x(0)=x_0.}
Therefore, the nominal stability assumption in Theorem~\ref{Tcir} means that  $A+B\Del_0(I-D\Del_0)^{-1}C$ is Hurwitz.
One often encounters the choice $\Del_0=0$ and the hypothesis that $A$ is Hurwitz.

To clarify the naming of the result, let $\nuu=\ny=1$, which implies $P\in\S^2$. Then, the set
$$
{\Bbb D}:=\left\{\la\in\C\st
\mat{c}{\la\\1}^*\!\!P\mat{c}{\la\\1}<0\right\}
$$
is an open disk in the complex plane different from $\emptyset$ and $\C$ iff $\det(P)<0$. It
could possibly degenerate to a half-plane as a disk in $\C\cup\{\infty\}$ \cite[Lemma A.1]{Sch05}.
Since \r{fdicir} translates into $\sy(\io)\in{\Bbb D}$ for all $\om\in\R\cup\{\infty\}$, the FDI
means that the Nyquist curve of $\sy$ is located in ${\Bbb D}$.

Although the article could focus on addressing the relation of this formulation of the circle criterion with so many others in the literature \cite{JayLog11}, only a few aspects related to dissipativity and integral quadratic constraints are mentioned in ``\nameref{sb:pcir}.''

\section[Dissipativity]{Dissipativity}\label{Sdiss}

Many excellent articles and monographs about dissipativity are available  \cite{HadChe08,ArcMei16,Van17,BroLoz19}.
The goal of this chapter is to collect some essential ingredients of dissipativity theory for robustness analysis as developed in this article, with a special emphasis on covering strict and nonstrict dissipativity in parallel.

\subsection{Nonlinear systems and supply rates}

This section recapitulates key dissipativity concepts for nonlinear systems.
With a vector field $f:\R^\nx\times \R^\nuu\to\R^\nx$ and a read-out map $g:\R^\nx\times \R^\nuu\to \R^\ny$,  consider
\eql{nlsy}{\dot x=f(x,u),\ \ y=g(x,u).}
An {\em admissible trajectory} of \r{nlsy} is a triple
$(x,u,y)\in\Se^\nx\times \Se^\nuu\times\Se^\ny$ that satisfies \r{nlsy}, again with the convention that the state trajectory $x$
is continuous and has a right derivative $\dot x$ in $\Se^n$.
Conditions that guarantee the existence of a solution of
the differential equation for some input $u\in\Se^\nuu$, and its uniqueness
if imposing $x(0)=x_0\in\R^\nx$, can be found in textbooks such as
\cite{LogRya14}. It causes no harm for our developments to assume that $f$ and $g$ are defined globally.

This brings us to the key notion of dissipativity, as introduced by Jan Willems in his celebrated work \cite{Wil72a,Wil72b,Wil71}, accompanied by the notion of strict dissipativity.

\definition{The system \r{nlsy} is {\em dissipative} with respect to the {\em supply rate} $s:\R^\nuu\times \R^\ny\to\R$ if there exists a {\em storage function} $V:\R^\nx\to\R$ such that the {\em dissipation inequality}
\equ{\tag{DI}\label{di}
V(x(t_2))\leq V(x(t_1))+\int_{t_1}^{t_2} s(u(t),y(t))\,dt
}
holds for all admissible trajectories and all $0\leq t_1\leq t_2$. Moreover, \r{nlsy} is {\em strictly dissipative} with respect to $s$ if there exist $V:\R^\nx\to\R$ and $\ve>0$ such that the
{\em strict dissipation inequality}
\equ{\tag{SDI}\label{sdi}    V(x(t_2))\ \leq \  V(x(t_1))  + \int_{t_1}^{t_2}s(u(t),y(t))\,dt
     -\varepsilon \int_{t_1}^{t_2}\|x(t)\|^2+\|u(t)\|^2\,dt
}
is satisfied for all admissible trajectories and all time instants $0\leq t_1\leq t_2$.
}
Since the system description and the supply rate do not explicitly depend upon time, it causes no loss of generality if taking $t_1=0\leq T=t_2$ in this definition.
Notably, both dissipation inequalities impose convex constraints on the pair $(V,s)$, which turns out to be essential for the construction of computational robust stability tests, as seen later.

Dissipativity is interpreted as a generalization of the law of conservation of energy.
The energy $V(x(t_1))$ stored in the system at time $t_1$ plus the supplied energy $\int_{t_1}^{t_2} s(u(t),y(t))\,dt$
in the time-interval $[t_1,t_2]$ is not smaller than the stored energy $V(x(t_2))$ at a later time $t_2$.
Strict dissipativity guarantees a strict decrease of the energy by an amount that is proportional to the
energy of the state trajectory plus that of the input signal on the interval $[t_1,t_2]$.
In contrast to Willems' original definition, it is not required that $V$ is bounded from below \cite{WilTak07}.

The global trajectory-based definitions of dissipativity can be replaced by static and local differential versions if $s$ is continuous
and $V$ is continuously differentiable
with the derivative $\partial V:=(\frac{\partial V}{\partial x_1},\ldots, \frac{\partial V}{\partial x_n})$. This  requires the following mild (standing) assumption on \r{nlsy}: Any pair $(x_0,u_0)\in\R^\nx\times \R^\nuu$ admits an admissible trajectory with $(x(0),u(0))=(x_0,u_0)$. The elementary proof of the following result can be found in \cite{Van17}, for example.

\theorem{\label{Tddis}
Let $V$ be continuously differentiable and $s$ be continuous. Then, the dissipation inequality \r{di} holds for all admissible trajectories iff the  differential dissipation inequality
\equ{\tag{DDI}\label{ddi}
\partial V(x)f(x,u)\leq s(u,g(x,u))\text{\ \ for all\ \ }(x,u)\in\R^\nx\times\R^\nuu
}
is satisfied. Similarly, the strict dissipation inequality \r{sdi} is equivalent to
\equ{\tag{SDDI}\label{sddi}
\partial V(x)f(x,u)\leq s(u,g(x,u))-\ve(\|x\|^2+\|u\|^2)\text{\ \ for all\ \ }(x,u)\in\R^\nx\times\R^\nuu.
}}

\skipthis{
\proo{We prove the second statement only, since the first is obtained by choosing $\ve=0$.
Take any $(x_0,u_0)\in\R^\nx\times \R^\nuu$ and a corresponding admissible trajectory
$(x,u,y)$ of \r{nlsy} with  $(x(0), u(0)) = (x_0, u_0)$. For all $t_1=0<t=t_2$, the strict dissipation inequality \r{sdi} implies
$$
\frac{V(x(t))-V(x(0))}{t}\leq \frac{1}{t}\int_{0}^{t} s(u(\tau),y(\tau))\,d\tau-\frac{\ve}{t}\int_0^t \|x(\tau)\|^2+\|u(\tau)\|^2\,d\tau.
$$
By taking the limit $t\searrow 0$ we infer
$$
\partial V(x(0))f(x(0),u(0))=\frac{d}{dt} V(x(t))|_{t=0}\leq s(u(0),y(0))-\ve(\|x(0)\|^2+\|u(0)\|^2),
$$
where we use the chain-rule on the left. Since $(x(0),u(0))=(x_0,u_0)$ and $y(0)=g(x_0,u_0)$, \r{sddi} follows because
$(x_0,u_0)$ was taken arbitrarily.

To prove the converse, let $(x,u,y)$ be any admissible system trajectory. We then infer from \r{sddi} with the system equations \r{nlsy}  that
$$
\partial V(x(t))\dot{x}(t)=\partial V(x(t))f(x(t),u(t))\leq s(u(t),y(t))-\ve(\|x(t)\|^2+\|u(t)\|^2)\text{\ \ for all\ \ }t\geq 0.
$$
Integration over $[t_1,t_2]$ with $0\leq t_1\leq t_2$ leads to \r{sdi}.\epro
}}

It is not exaggerated to claim that a multitude of concrete applications of dissipativity rests on the step of moving from \r{ddi} to \r{di} [or from \r{sddi} to \r{sdi}], just by integration on any time interval $[t_1,t_2]\subset[0,\infty)$ along a system trajectory. In fact, $s$ often encodes some desired trajectory-based property of the system \r{nlsy}, which can be expressed as a consequence of \r{di}, by possibly imposing extra conditions on $V$ or restricting the considered class of trajectories. Then, \r{ddi} is proposed as an algebraically verifiable condition to guarantee \r{di}.

\subsection{Example: Characterizing energy gain bounds and extensions}

For a most prominent example, let \r{nlsy} have the equilibrium
$(x_e,u_e)\in\R^n\times\R^m$, which means  that $f(x_e,u_e)=0$. Moreover, with $\ga>0$, choose
$s(u,y):=\ga^2\|u\|^2-\|y\|^2.$ Then, \r{ddi}   with  $0=t_1\leq t_2=T$ implies that
\eql{gdi}{
V(x(T))\leq V(x(0))+\int_{0}^{T} \ga^2\|u(t)\|^2-\|y(t)\|^2\,dt\te{for all}T\geq 0.
}
This is an often-seen energy-gain property for the nonlinear system that, by itself, can be read and interpreted in various
ways, two of which are highlighted.

If $V$ has a global minimum at $x_e$, we may replace $V$ with $V-V(x_e)$. This assures $V(x_e)=0$ and
 $V(x)\geq V(x_e)=0$ for all $x\in\R^\nx$. For all admissible trajectories,
\r{gdi} then implies
$$
\int_{0}^{T} \|y(t)\|^2\,dt\leq V(x(0))+\ga^2\int_{0}^{T} \|u(t)\|^2\,dt\te{for all}T\geq 0.
$$
For finite energy inputs $u\in\Se^\nuu$, this shows
$\int_{0}^{T} \|y(t)\|^2\,dt\leq V(x(0))+\ga^2\|u\|_2^2$ for all $T\geq 0$, which not only assures that $y$ has finite energy, but
also leads to
\eql{gain}{
\|y\|_2^2\leq V(x(0))+\ga^2\|u\|_2^2\te{for all}u\in\Se^\nuu\te{with}\|u\|_2<\infty.
}
In the terminology of \cite[p. 41]{DesVid75}, this means that the  $\Le_2$-gain of \r{nlsy} is bounded by $\ga$. The so-called offset
$V(x(0))$ is related to the system's initial condition $x(0)$ by the storage function. Imposing additional growth conditions on $V$ leads to inequalities that can be directly expressed in terms of $x(0)$ and $u$. It is also
often prescribed to initialize the system at the equilibrium as $x(0)=x_e$ such that $V(x(0))=V(x_e)=0$. Then, \r{gain} guarantees that the gain of
\r{nlsy} is not larger than $\ga$ according to \cite[p. 232]{Zam66}.

Exactly the same arguments show that strict $s$-dissipativity of \r{nlsy} leads to
\eql{sgain}{\ve\|x\|_2^2+\|y\|_2^2\leq V(x(0))+(\ga^2-\ve)\|u\|_2^2\te{for all}u\in\Se^\nuu.}
If $u$ is of finite energy, this implies the same for all state and output trajectories, and the energy gain of \r{nlsy} is strictly smaller than $\ga$.

Apart from the interpretation of the consequences, nothing depends on the choice of the supply rate. Following \cite{Wil72b,Wil71,HilMoy77}, a still limited but much larger flexibility is gained by taking a homogeneous quadratic function defined as
\eql{qsr}{
s_P(u,y)=\mat{c}{y\\u}^T\!\!P\mat{c}{y\\u}\te{with}P=\mat{cc}{Q&S\\S^T&R}\in\S^{\ny+\nuu}.
}
The most popular choices are $-P_{\gai}$ and $P_{\pas}$ with the ``bounded real'' and ``positive real'' matrices
\eql{P}{
P_{\gai}:=\mat{cc}{I_\ny&0\\0&-I_\nuu}\te{and}
P_{\pas}:=\mat{cc}{0&\frac{1}{2}I_\nuu\\\frac{1}{2}I_\nuu&0},
}
respectively, where  $\nuu=\ny$ in the latter case.
The dissipation inequality for $s_{P_{\pas}}$ reads as
$$
V(x(T))\leq V(x(0))+\int_{0}^{T} u(t)^Ty(t)\,dt\te{for all}T\geq 0.
$$
Again, varying assumptions on $V$ and the system initialization lead to the different
conclusions around the classical notion of passivity in the literature. This special case has been addressed in much detail in \cite{Wil72b}. Other choices of $P$ capture the concepts of strict input- and output-passivity and many others \cite{HilMoy77}. The latter article also addresses consequences of the fact that the nonlinear system is input-affine.

\subsection{Linear systems and quadratic supply rates}

Dissipativity for  linear systems and quadratic supply rates is of particular importance for this article. This is the reason why we recapitulate the key results in much detail and under weak assumptions.
Apart from the KYP Lemma, all statements are proven with complete arguments.

The choice of linear maps $f$ and $g$ in \r{nlsy} leads to a linear control system
\eql{lsy}{\dot x=Ax+Bu,\ \ y=Cx+Du.}
For a quadratic supply rate $s_P$ as in \r{qsr}, it is natural to pick  quadratic storage functions
\eql{qsf}{
V(x)=x^T\Xs x\te{for}x\in\R^\nx\te{with a describing matrix}\Xs\in\S^\nx
}
to characterize dissipativity. Trivially, (strict) dissipativity with a quadratic storage function implies (strict) dissipativity in general. For a linear system and a quadratic supply rate, the converse is true as well.
However, for nonstrict dissipativity, it is instrumental to assume that \r{lsy} is controllable, while strict dissipativity does not require any such hypothesis.

\theorem{\label{Tqsf}
If the linear system \r{lsy} is strictly $s_P$-dissipative with respect to the quadratic supply rate \r{qsr}, then it is strictly $s_P$-dissipative with a homogeneous quadratic storage function \r{qsf}.
For controllable systems \r{lsy}, this also holds for (nonstrict) $s_P$-dissipativity.}

This is a consequence of the celebrated Kalman-Yakubovich-Popov lemma (Theorem~\ref{Tkyp}).

For quadratic storage functions \r{qsf}, it has been emphasized in \cite{Wil71} and is easy to establish that the  \r{ddi}
translates into the following LMI in the defining matrix $\Xs$:
\eql{lmidi}{\mat{cc}{A&B\\I&0}^T\!\!\mat{cc}{0&\Xs\\\Xs&0}\mat{cc}{A&B\\I&0}\cle
\mat{cc}{C&D\\0&I}^T\!\!\!P\mat{cc}{C&D\\0&I}.}
Similarly, the \r{sddi} leads to the strict LMI
\eql{lmisdi}{\mat{cc}{A&B\\I&0}^T\!\!\mat{cc}{0&\Xs\\\Xs&0}\mat{cc}{A&B\\I&0}\cl
\mat{cc}{C&D\\0&I}^T\!\!\!P\mat{cc}{C&D\\0&I}.}

\theorem{\label{Tlmi}The quadratic storage function $V(x)=x^T\Xs x$ with $\Xs\in\S^\nx$ satisfies the \r{ddi}  for the linear system \r{lsy} and the supply rate $s_P$
iff $\Xs$ satisfies the LMI \r{lmidi}. The same holds for
the strict dissipation inequality \r{sddi} and the strict LMI \r{lmisdi}.}

\proo{The proof requires to observe that the \r{ddi} or \r{sddi}  involve inequalities for homogeneous quadratic functions, which
equivalently translate into the corresponding matrix inequality constraints for their Hessians.
Explicitly,  for $V$ in \r{qsf},
we infer $\partial V(x)f(x,u)=2x^T\Xs(Ax+Bu)$.
In view of \r{qsr}, the \r{ddi} is
$$
\mat{cc}{Ax+Bu\\x}^T\!\!\mat{cc}{0&\Xs\\\Xs&0}\mat{cc}{Ax+Bu\\x}
\leq\mat{cc}{Cx+Du\\u}^T\!\!P\mat{cc}{Cx+Du\\u}
$$
for all $x\in\R^\nx$ and $u\in\R^\nuu$. This is equivalent to \r{lmidi}.
In the same way, the the \r{sddi} translates into the strict LMI
\r{lmisdi}.
}

Frequency domain inequalities (FDIs) emerge if driving  \r{lsy} with periodic complex input functions
$u(t)=u_0e^{\io t}$ for $t\geq 0$, with frequency $\om\in\R$ and where $u_0\in\C^\nuu$ is a complex amplitude vector. Such an input generates
the state and output responses $x(t)=x_0e^{\io t}$ and $y(t)=y_0e^{\io t}$ for $t\geq 0$, if $x_0\in\C^\nx$ and $y_0\in\C^\ny$
satisfy the algebraic equations
$(A-\io I)x_0+Bu_0=0$, $y_0=Cx_0+Du_0$.
This motivates the introduction of the subspace
$$
\c{\sy}(\la):=\{(x_0,u_0)\in\C^\nx\times\C^\nuu\st (A-\la I)x_0+Bu_0=0\}\te{for}\la\in\C.
$$
For strict dissipativity, also consider $\c{\sy}(\infty):=\{(x_0,u_0)\in\C^{\nx+\nuu}\st x_0=0\}.$

\lemma{\label{Lfdi}
If \r{lsy} is $s_P$-dissipative, then
\eql{fdi}{
0\leq s_P(u_0,Cx_0+Du_0)
\te{for all}(x_0,u_0)\in\c{\sy}(\io)\te{and all}\om\in\R.
}
If \r{lsy} is strictly  $s_P$-dissipativity, then
\eql{sfdi}{
0<s_P(u_0,Cx_0+Du_0)
\te{for all}(x_0,u_0)\in\c{\sy}(\io)\setminus\{0\}\te{and all}\om\in\R\cup\{\infty\}.
}
}

Before entering the proof, we address how these conditions are verifiable.
For this purpose, choose a basis matrix
$V_\la\in\C^{(\nx+\nuu)\times p}$ of the subspace $\c{\sy}(\la)$.
For $\la=\infty$, take $V_\infty=\col(0,I_\nuu)$. 
If $\la$ is not an eigenvalue of $A$, note that $(A-\la I)x_0+Bu_0=0$ is equivalent to
$x_0=(\la I-A)^{-1}Bu_0$. Then, $\c{\sy}(\la)$ has the basis
\eql{bas}{
V_\la:=\mat{cc}{(\la I-A)^{-1}B\\I_\nuu}\te{for}\la\in\C\cup\{\infty\},\ \la\not\in\eig(A).
}
Here, $\eig(M)$ denotes the set of complex eigenvalues of some real or complex matrix $M$.
This implies  that \r{sfdi} just requires the verification of
$$
0\cl V_{\io}^*\mat{cc}{C&D\\0&I}^T\!\!P\mat{cc}{C&D\\0&I}V_{\io}\te{for all}\om\in\Ri,
$$
while \r{fdi} translates into the nonstrict inequality for $\om\in\R$.

\proo{
Let \r{lsy} be $s_P$-dissipative and choose $\om\in\R$ as well as $(x_0,u_0)\in\c{\sy}(\io)$.
With $y_0=Cx_0+Du_0$, the signals
$u(t)=u_0e^{\io t}$, $x(t)=x_0e^{\io t}$, and $y(t)=y_0e^{\io t}$ for $t\geq 0$
define a complex trajectory of \r{lsy}. Their real and imaginary parts give two real trajectories.

For $\om>0$, take $0=t_1<t_2=\frac{2\pi}{\om}=:T$ to obtain $x(T)=x(0)$.
Therefore, \r{di} implies
$\int_0^{T} s_P(\Re(u(t)),\Re(y(t))\,dt\geq 0$ and $\int_0^{T} s_P(\Im(u(t)),\Im(y(t))\,dt\geq 0$.
Since $P$ is real symmetric, $z\in\C^{\ny+\nuu}$ gives
$z^*Pz=(\Re(z)^T-i\Im(z)^T)P(\Re(z)+i\Im(z))=
\Re(z)^TP\Re(z)+\Im(z)^TP\Im(z)$. Therefore, the latter two inequalities can be added to obtain
$$0\leq
\frac{1}{T}
\int_0^T \mat{c}{y_0e^{i\om t}\\u_0e^{i\om t}}^*\!\!P\mat{c}{y_0e^{i\om t}\\u_0e^{i\om t}}\,dt=
\frac{1}{T}
\int_0^T \mat{c}{y_0\\u_0}^*\!\!P\mat{c}{y_0\\u_0}\,dt=s_P(u_0,y_0),
$$
where we exploited $(e^{i\om t})^*e^{i\om t}=1$.
For $\om=0$, the trajectories are constant and the latter inequality is obtained with any $T>0$. For $\om<0$, we take instead $0=t_1<t_2=-\frac{2\pi}{\om}=:T$ to arrive at the same conclusion. This completes the proof of the first claim.

If \r{lsy} is strictly $s_P$-dissipative, identical arguments lead to
$\ve\|x_0\|^2+\ve\|u_0\|^2\leq s_P(u_0,y_0),$
which shows $0<s_P(u_0,y_0)$ for $(x_0,u_0)\neq 0$.
The case $\om=\infty$ is handled with a limiting argument. Assume $0\neq (0,u_0)\in \c{\sy}(\infty)$.
Then, take a sequence of finite frequencies with $\om_\nu\to\infty$ for $\nu\to\infty$ and such that $\io_\nu$ is no eigenvalue of $A$. With $x_\nu:=(i\om_\nu I-A)^{-1}Bu_0$, one infers $(x_\nu,u_0)\in\c{\sy}(\io_\nu)$ and $(x_\nu,u_0)\to (0,u_0)$ for
$\nu\to\infty$. Since $\ve\|x_\nu\|^2+\ve\|u_0\|^2\leq s_P(u_0,Cx_\nu+Du_0)$, we can take the limit $\nu\to\infty$ to conclude
$0<s_P(u_0,y_0)$ because of $u_0\neq 0$. This proves the
second claim.
}

The celebrated Kalman-Yakobovich-Popov (KYP) lemma now completes the cycle of arguments by linking the FDIs \r{fdi}-\r{sfdi}
to the dissipation LMIs \r{lmidi}-\r{lmisdi}, respectively. Although this result can be traced back to the 1960s \cite{Yak62}, we rely on a
version in \cite{BalVan03} that has been proven by duality in convex optimization.

\theorem{\label{Tkyp}The strict dissipation LMI \r{lmisdi} has a symmetric solution iff \r{sfdi} holds.
If the system \r{lsy} is controllable, the nonstrict dissipation LMI \r{lmidi} is feasible iff \r{fdi} holds.}

\proo{Apply  Lemmas 1 and 2 in \cite{BalVan03} for the strict the
nonstrict versions, respectively.
}

The ``if'' statements in this KYP Lemma are traditionally considered to be the difficult parts in all proofs.
Note that the ``only if'' directions are direct consequences of a combination of Theorem~\ref{Tlmi} with Lemma~\ref{Lfdi}.
Moreover, it also follows that the dissipativity characterizations in Lemma~\ref{Lfdi} are exact.

Theorem \ref{Tqsf} can be proved using the KYP Lemma. To this end, let \r{lsy} be $s_P$-dissipative with a
general storage function. By Lemma~\ref{Lfdi}, infer that the FDI \r{fdi} holds true. By the KYP Lemma (Theorem \ref{Tkyp}), the LMI \r{lmidi} has a symmetric solution $\Xs$. By Theorem~\ref{Tlmi}, this solution defines a quadratic storage function \r{qsf} that
certifies $s_P$-dissipativity of \r{lsy}. Exactly the same arguments work for strict dissipativity.

It  is remarkable that the results for strict dissipativity do not require any assumptions on the data matrices describing the system \r{lsy} or the supply rate \r{qsr}. For dissipativity, controllability of $(A,B)$ in Theorem \ref{Tkyp} cannot be weakened without additional assumptions. If $s_P(I,D)\cg 0$, the LMI turns into an algebraic Riccati inequality,
the theory of which is understood under much weaker assumptions on $(A,B)$ (see \cite{Sch95a}, for example).
Note that contributions that address variations of Theorem~\ref{Tkyp} abound in the literature.

With the transfer matrix $\sy(s)=C(sI-A)^{-1}B+D$ of the controllable linear system \r{lsy}, the FDI \r{fdi} is equivalent to the more traditional version formulated as
\eql{fdidi}{
0\cle \mat{c}{\sy(\io)\\I_\nuu}^*\!\!P\mat{c}{\sy(\io)\\I_\nuu}
\te{for all $\om\in \R$ such that $\io$ is no pole of}\sy.}
If $A$ has no eigenvalues in $i\R$, then \r{sfdi} is equivalent to
\eql{fdisdi}{
0\cl \mat{c}{\sy(\io)\\I_\nuu}^*\!\!P\mat{c}{\sy(\io)\\I_\nuu}
\te{for all $\om\in \Ri$.}
}
The latter follows from the fact that $V_{\io}$ is a basis of $\c{\sy}(\io)$ for all $\io\in i\R\cup\{\infty\}$ by \r{bas}.
\lemma{\label{LfdiII} If $(A,B)$ is controllable, then \r{fdi} $\!\Leftrightarrow\!$ \r{fdidi}.
If $\eig(A)\cap i\R\!=\!\emptyset$,
then
\r{sfdi} $\!\Leftrightarrow\!$ \r{fdisdi}.
}
\proo{
To show \r{fdi} $\Rightarrow$ \r{fdidi}, select any $\io_0\in i\R$ that is not a pole of $\sy$. Then, choose a sequence $\om_\nu$ with $\om_\nu\to\om_0$ for $\nu\to\infty$ such that $i\om_\nu\not\in\eig(A)$. Since $V_{\io_\nu}$ in \r{bas} is a basis of $\c{\sy}(\io_\nu)$, \r{fdi} implies
$0\cle s_P(I_\nuu,\sy(\io_\nu))$ for all $\nu\in\N$. Due to $\sy(i\om_\nu)\to \sy(\io_0)$ for $\nu\to\infty$,
the limit can be taken to infer $0\cle s_P(I_\nuu,\sy(\io_0))$ and, hence, \r{fdidi}.

To show \r{fdi} $\Leftarrow$ \r{fdidi}, suppose that $(A,B)$ is controllable. Take any $\io_0\in i\R$ and $(x_0,u_0)\in\c{\sy}(\io_0)$. Since $(A-\io_0 I\ B)$ has full row rank, we can find real matrices $(\t C\ \t D)$ such that
$$
S_{\om}:=\mat{cc}{A-\io I&B\\\t C&\t D}\te{is invertible for}\om=\om_0.
$$
Choose a sequence $\om_\nu\neq \om_0$ with $\om_\nu\to\om_0$ for $\nu\to\infty$ such that $S_{\om_\nu}$ is invertible
and $i\om_\nu\not\in\eig(A)$. Then, define
$$
\mat{c}{x_\nu\\u_\nu}:=S_{\om_\nu}^{-1}S_{\om_0}\mat{c}{x_0\\u_0}.
$$
Since $S_{\om_\nu}\to S_{\om_0}$, we conclude $\col(x_\nu,u_\nu)\to \col(x_0,u_0)$ for $\nu\to \infty$.
Moreover, note that $(A-\io_\nu I)x_\nu+Bu_\nu=0$ for all $\nu\in\N$.
Since $i\om_\nu\not\in\eig(A)$, infer $x_\nu=(\io_\nu I-A)^{-1}Bu_\nu$ and thus $Cx_\nu+Du_\nu=\sy(\io_\nu)u_\nu$ for all $\nu\in\N$.
Then, \r{fdidi} shows $0\leq s_P(u_\nu,\sy(\io_\nu)u_\nu)$ and thus $0\leq s_P(u_\nu,Cx_\nu+D u_\nu)$ for all $\nu\in\N$.
Taking the limit $\nu\to\infty$ gives $0\leq s_P(u_0,Cx_0+Du_0)$.  Since $\io_0\in i\R$ and $(x_0,u_0)\in\c{\sy}(\io_0)$ were taken arbitrarily, the proof is completed.
}

\subsection{Sign-constrained storage functions}

So far, all given results address the existence of indefinite storage functions. For strict dissipativity, the existence of a positive definite storage functions is guaranteed by suitable stability assumptions as follows.

\theorem{\label{Tlmipos}
Suppose that there exist some $\Delta_0\in\R^{\nw\times\nz}$ with $s_P(\Delta_0,I)\cle 0$, and let $X=X^T$ satisfy \r{lmisdi}. Then, $I-D\Del_0$ is invertible. Moreover,  all solutions of the strict dissipation LMI \r{lmisdi} are positive definite iff $A+B\Del_0(I-D\Del_0)^{-1}C$ is Hurwitz.
}

\proo{To show $\det(I-D\Del_0)\neq 0$, note that \r{lmisdi} implies $s_P(I,D)\cg 0$. Therefore,
\eql{h00}{
s_{P}(\Del_0z,z)\leq 0\te{for all}z\in\R^\nw
\te{and}
s_{P}(w,Dw)>0\te{for all}w\in\R^{\nw}\setminus\{0\}.
}
Hence, $(I-D\Delta_0)z=0$ for $z\neq 0$ implies
with $w:=\Delta_0z$ that $z=Dw$ and $w\neq 0$; then \r{h00} leads to
$s_{P}(w,z)\leq 0$ and $s_{P}(w,z)>0$, a contradiction. 

To show the claimed equivalence, observe that
$\t C:=(I-D\Del_0)^{-1}C$ satisfies $C+D\Del_0\t C=\t C$.
Right-multiplying \r{lmisdi} with $\col(I,\Del_0\t C)$ and left-multiplying the transpose gives
\eqn{
(A+B\Del_0\t C)^T\Xs+\Xs(A+B\Del_0\t C)\cl \t C^Ts_P(\Delta_0,I)\t C,
}
and thus $(A+B\Del_0\t C)^T\Xs+\Xs(A+B\Del_0\t C)\cl 0$.
This inequality completes the proof, since it shows  $\Xs\cg 0$ iff
$A+B\Del_0\t C=A+B\Del_0(I-D\Del_0)^{-1}C$ is Hurwitz.
}

The formulation of Theorem~\ref{Tlmipos} simplifies if we can choose $\Delta_0=0$. Then, $s_P(0,I)\cle 0$ means that the left-upper block of $P$ is negative semidefinite.
The given proof also covers the following consequence for the nonstrict dissipation LMI.


\corollary{\label{Clmipos} Suppose that the left-upper block of $P$ is negative semidefinite. If $A$ is Hurwitz then all solutions of the nonstrict dissipation LMI \r{lmidi} are positive semidefinite.
Moreover, all solutions of the strict dissipation LMI \r{lmisdi} are positive definite iff $A$ is Hurwitz.
}

A concise summary of the insights in this section for the supply rate $-s_P$ is given in ``\nameref{sb:dis}.'' The change of sign leads to formulations that are compatible with the subsequent discussion on robustness.

\subsection{Ramifications}

Everything discussed up to now can also be established for systems defined on the discrete time set $\{0,1,2,3,\ldots\}$ with hardly any modifications. Next to changing the time set, the most essential adaptation is to replace the imaginary axis $i\R$ with the unit circle $\{z\in\C\st |z|=1\}$ for frequency domain inequalities. The corresponding LMIs are then obtained
by
$$\text{replacing\ \ }\mat{cc}{0&\Xs\\\Xs&0}\te{with}\mat{cc}{\Xs&0\\0&-\Xs}.$$
This is one reason for the given formulations of the dissipation LMIs without working out the matrix products explicitly,
as is the habit in so many other texts. This version nicely exhibits the system matrices, the supply rate, and the
storage function matrix in a separated fashion, which helps to trace these ingredients in more complicated cases (as developed in the sequel).

Various further extension can be routinely established. For example, still on the time set $[0,\infty)$,  one
can consider time-varying systems
$\dot x=f(t,x,u)$, $y=g(t,x,u)$ and study dissipativity for time-varying supply rates $s(t,u,y)$. As the only essential modification, one then works with nonautonomous storage functions $V(t,x)$. Because of $\frac{d}{dt}V(t,x(t))=
\partial_tV(t,x(t))+\partial_xV(t,x(t))\dot{x}(t)$, this leads to the
differential dissipation inequality
$$
\partial_tV(t,x)+\partial_xV(t,x)f(t,x,u)+s(t,u,g(t,x,u))\leq 0\te{for all}t\geq 0,\ (x,u)\in\R^n\times\R^\nuu.
$$
For time-varying linear systems
$\dot x=A(t)x+B(t)u$, $y=C(t)x+D(t)u$,
one often uses quadratic supply rates $s_{P(t)}(u,y)$ with a time-varying matrix $P(t)$. With a quadratic storage function $V(t,x):=x^T X(t)x$,
the DDI then turns into a differential LMI, which means that $X(.)$ satisfies
$$\arraycolsep.4ex
\mat{cc}{A(t)&B(t)\\I&0}^T\!\!\mat{cc}{0&X(t)\\X(t)&\dot{X}(t)}\mat{cc}{A(t)&B(t)\\I&0}
\cle\mat{cc}{C(t)&D(t)\\0&I}^T\!\!P(t)\mat{cc}{C(t)&D(t)\\0&I}\text{\ for\ }t\geq 0.
$$
Strict dissipativity results (and variations thereof) can be established analogously.
For linear time-varying systems (for example),
strict differential LMIs (characterizing strict dissipativity) are exploited in \cite{Sch96} for multiobjective controller synthesis. 

\section[Robustness analysis with dynamic IQCs]{Robustness analysis with dynamic integral
quadratic constraints }\label{RS}

Building on the ideas in ``\nameref{sb:pcir},'' the current section
covers three substantial extensions. First, we move from
simple static nonlinear mappings in the feedback loop to inclusions that comprise
classes of single- or multivalued general static or dynamic uncertainties.
Following \cite{MoyHil78,SafAth81}, the properties of the uncertainties are captured
by families of multipliers that define valid so-called static IQCs.
The incorporation of dynamics in the IQCs substantially widens the scope of and strengthens dissipativity-based robustness tests \cite{MegRan97}. In contrast to \cite{MegRan97}, we follow \cite{SchVee18} and work with the notion of IQCs with terminal cost. This leads to a seamless link of dissipation theory and IQC theory, as exposed in detail in ``\nameref{sb:iqc}.'' Finally, this section covers the characterization of robust performance with a rather detailed exposition of the wealth of possible interpretations.

\subsection[Robust stability and static IQCs]{Robust stability and static
integral
quadratic constraints }
\label{Srsh}

Consider the feedback loop
\eql{fbg}{\dot x=Ax+Bw,\ \ z=Cx+Dw,\ \ w\in\Des(z),\te{and}x(0)=x_0}
as depicted on the left in Fig.~\ref{figrp}.
If compared to \r{fb}, the linear system with transfer matrix $\sy$ remains unchanged, while
the uncertainty $\Des$ can now be any multivalued mapping that
assigns a whole set of signals $\Des(z)\subset\Se^\nw$ to
a given  $z\in\Se^\nz$.
Still, $(x,w,z)\in\Se^\nx\times\Se^\nw\times\Se^\nz$ is a response of the loop if the signals satisfy
\r{fbg}, and the loop is stable if there exists some $\ga>0$ such that \r{sta} holds true along all responses.
It is assumed that
\eql{nom}{0\in\Des(z)\te{for all}z\in\Se^\nz.}
Then, the responses of \r{fbg} comprise those of the nominal system $\dot x=Ax$, $z=Cx$, $w=0$, and $x(0)=x_0$.
Therefore, nominal stability of the loop is equivalent to $A$ being Hurwitz.

As an example, the circle criterion involves the uncertainty
$$\Des(z):=\{w\in\Se^\nw\st w=\vp(z),\ \vp:\R^\nz\to\R^\nw \text{ satisfies \r{qccir}}\} \te{for}z\in\Se^\nz.$$
This indicates why it is natural to work with set-valued uncertainties,
even for feedback loops \r{fb} that are defined with  single-valued mappings \cite{Zam66}.
If the zero map $\vp=0$ satisfies \r{qccir}, then $0\in\Des(z)$ for all $z\in\Se^\nz$,
which illustrates how this property emerges in practice.

As discussed in ``\nameref{sb:pcir},'' absolute stability tests exploit the fact that all signals $(w,z)\in\Se^\nw\times\Se^\nz$ with $w\in\Des(z)$ satisfy a so-called hard static IQC
\eql{hiqc}{\int_{0}^{T}s_P(w(t),z(t)) \, dt\geq 0\te{for all}T\geq 0}
for the multiplier $P\in\S^{\nz+\nw}$. ``Hard'' refers to the fact that the inequality holds on the finite time interval $[0,T]$ for all $T\geq 0$, while soft IQCs require it to hold on $[0,\infty)$. ``Static'' refers to the fact that the supply rate
$s_P(w,z)$ is a nondynamic function of the signal $(w,z)$, while dynamic IQCs involve filtered versions thereof.
Observe that \r{nom} and \r{hiqc} imply \r{iqcnom} to hold for $\Delta_0=0$.

As addressed for the circle criterion, the generation of powerful robustness tests relies on the construction of a rich family $\bm{P}\subset\S^{\nz+\nw}$ of multipliers with an LMI representation such that $\Des$ satisfies a hard IQC for any $P\in\bm{P}$.
Then, the proof of the circle criterion needs no adaptation to cover the following more general result.
\corollary{\label{Crsh}
The loop  \r{fbg} is stable if there exist $P\in\bm{P}$ and $\Xs\in\S^n$ which satisfy
\eql{rslmi}{
\Xs\cg 0\te{and}
\mat{cc}{A&B\\I&0}^T\!\!\mat{cc}{0&\Xs\\\Xs&0}\mat{cc}{A&B\\I&0}+\mat{cc}{C&D\\0&I}^T\!\!P\mat{cc}{C&D\\0&I}\cl 0.
}
Moreover, for fixed $P\in\bm{P}$, \r{rslmi} is feasible in $X$ iff $A$ is Hurwitz and the FDI \r{fdicir} holds.
}

\subsection{Robust stability and dynamic
integral
quadratic constraints }

Dynamic IQCs correspond to working with a multiplier family $\Psi^*\bm{P}\Psi$ defined by a set of symmetric matrices $\bm{P}\subset\S^p$ and some fixed filter
that is described by a transfer matrix $\Psi=\mat{cc}{\Psi_1&\Psi_2}$  of dimension $p\times (\nz+\nw)$.
With a state-space realization of $\Psi$, the response $v=\Psi_1z+\Psi_2w$ is defined through
\eql{psi}{
\left.\arr{rcl}{\dot \xi&=&A_\Psi \xi+B_{\Psi_1}z+B_{\Psi_2}w,\ \ \xi(0)=0,\\v&=&C_\Psi \xi+D_{\Psi_1}z+D_{\Psi_2}w.}\right.
}

\definition{\label{DiqcZ}The uncertainty $\Des$ satisfies an IQC
with terminal cost for the multiplier class $\Psi^*\bm{P}\Psi$ if for every $P\in\bm{P}$, there exists a terminal cost matrix $Z=Z^T$ such that
\eql{iqcZ}{
\int_0^T v(t)^T\mi v(t)\,dt-\xi(T)^TZ\xi(T)\geq 0\te{ for all}T\geq0
}
along all filtered signals $v=\Psi_1z+\Psi_2w$ with $(w,z)\in\Se^\nw\times\Se^\nz$ and $w\in\Des(z)$.
}

This notion was introduced in \cite{SchVee18} (with a sign-change of the terminal cost term for reasons of compatibility with the current article). It builds on the work in \cite{Sei15} for a fixed
multiplier $\Psi^*P\Psi$ with $P\in\bm{P}$, which formulates conditions under which an alternative factorization $\hat\Psi^*\hat P\hat\Psi$ of the given multiplier leads to an IQC with terminal cost zero. This translates into the fact that
the uncertainty satisfies an IQC with a terminal cost matrix $Z$, which is uniquely
determined by the stabilizing solution of the algebraic Riccati equation in Theorem~\ref{Tiqc2} defined by $\Psi$ and $P$.

In contrast and aligned with dissipativity theory, Definition~\ref{DiqcZ} allows
for a whole set of terminal cost matrices $Z$, even if $P$ is fixed. Moreover, the set of pairs $(P,Z)$ with
$P=P^T$ and $Z=Z^T$ that satisfy \r{iqcZ} is clearly convex, which is of key relevance for the generation of computational robustness tests. Still, this concept differs from a standard dissipativity constraint, since it only involves a partial storage function of the state of the filter defined by $Z$ and does not rely on any information about ``internal properties'' of $\Des$ [such as their state trajectories if described by a family of ODEs].
If $\Psi$ is the trivial filter $\Psi=I_{\nz+\nw}$, then \r{iqcZ} reduces to the hard static IQC \r{hiqc}.
If $\Psi$ is not trivial, the IQC is said to be dynamic, since the multipliers in $\Psi^*\bm{P}\Psi$ are frequency dependent and
define a dynamic supply rate.

The subsequent robust stability test involves the FDI
\eql{fdis}{
\mat{c}{\sy\\I}^*\!\!\Psi^*\mi\Psi\mat{c}{\sy\\I}=
\left[\Psi_1\sy+\Psi_2\right]^*\!\!\mi\left[\Psi_1\sy+\Psi_2\right]
\cl0
}
to hold for some $P\in\bm{P}$. According to Corollary~\ref{Csdi}, this FDI  relates to a strict dissipativity condition for the system $v=(\Psi_1G+\Psi_2)w$ and a supply rate in terms of $P$.
A natural state-space description of this system as the series interconnection of the filter $\Psi$ and $\col(\sy,I)$
is given by
\eql{syf}{
\arr{rcl}{
\mat{c}{\dot \xi\\\dot x}&=&
\mat{cc}{
A_\Psi&B_{\Psi_1}C\\
0     &A         }
\mat{c}{\xi\\x}+
\mat{cc}{
B_{\Psi_1}D+B_{\Psi_2}\\
B}w,\ \ \xi(0)=0,\\
v&=&
\mat{cc}{C_\Psi&D_{\Psi_1}C}
\mat{c}{\xi\\x}+
(D_{\Psi_1}D+D_{\Psi_2})w.
}
}
Therefore, the corresponding dissipation LMI reads as
\mul{\label{lmif}
\mat{cc|cc}{A_\Psi&B_{\Psi_1}C&B_{\Psi_1}D+B_{\Psi_2}\\0&A&B\hl I&0&0\\0&I&0}^T\!\!
\mat{cc|cc}{0&0&X_{1}&X_{12}\\0&0&X_{12}^T&X_{2}\hl X_{1}&X_{12}&0&0\\X_{12}^T&X_{2}&0&0}
\mat{cc|cc}{A_\Psi&B_{\Psi_1}C&B_{\Psi_1}D+B_{\Psi_2}\\0&A&B\hl I&0&0\\0&I&0}+\\
+\mat{cc|c}{C_\Psi&D_{\Psi_1}C&D_{\Psi_1}D+D_{\Psi_2}}^T\!\!P\mat{cc|c}{C_\Psi&D_{\Psi_1}C&D_{\Psi_1}D+D_{\Psi_2}}\cl 0.}
This leads to the following robust stability result, a variant of which appeared in \cite{SchVee18}.

\theorem{\label{Trs}
Let $\Des$ satisfy a finite-horizon IQC with terminal cost for the multiplier class $\Psi^*\bm{P}\Psi$. Then, \r{fbg} is robustly stable if there exists some $P\in\bm{P}$, a corresponding terminal cost matrix $Z=Z^T$, and some $X=X^T$ with \r{lmif} that are coupled as
\eql{pos}{
\mat{cc}{X_{1}+Z&X_{12}\\X_{12}^T&X_{2}}\cg 0.
}
}

\proo{Select any trajectory of \r{fbg}. By Corollary~\ref{Csdi} and \r{lmif}, there exists some $\ve>0$ such that
the given loop trajectory filtered by \r{psi} satisfies the dissipation inequality
\mun{
\mat{c}{\xi(T)\\x(T)}^T\!\!\mat{cc}{X_{1}&X_{12}\\X_{12}^T&X_{2}}\mat{c}{\xi(T)\\x(T)}+\int_0^Tv(t)^TPv(t)\,dt+\\+\ve\int_0^T\|\xi(t)\|^2+\|x(t)\|^2+\|w(t)\|^2\,dt\leq
\mat{c}{0\\x(0)}^T\!\!\mat{cc}{X_{1}&X_{12}\\X_{12}^T&X_{2}}\mat{c}{0\\x(0)}\te{for all}T\geq 0.
}
Since $w\in\Des(z)$ (and by the assumption on $\Des$), we can exploit \r{iqcZ} to bound
$\int_0^Tv(t)^T\mi v(t)\,dt$ from below by $\xi(T)^T Z\xi(T)$ for all $T\geq 0$.
This implies
\eqn{
\mat{c}{\xi(T)\\x(T)}^T\!\!\mat{cc}{X_{1}+Z&X_{12}\\X_{12}^T&X_{2}}\mat{c}{\xi(T)\\x(T)}+\ve\!\int_0^T\|x(t)\|^2+\|w(t)\|^2\,dt\leq
x(0)^TX_{2}x(0)
}
for all $T\geq 0$, which proves \r{sta} due to \r{pos} and completes the proof.
}

In contrast to other dissipativity results, Theorem~\ref{Trs} features a coupled positivity constraint on  $Z$ and $X$,
the terminal cost matrix of the uncertainty IQC \r{iqcZ} (viewed as a partial storage function for the filter's state)
and the storage function certifying dissipativity of the filtered linear system \r{syf}.
As illustrated in the middle Fig.~\ref{figrp}, the coupling is motivated by the trajectories of
both the uncertainty and of the linear system jointly driving the  filter dynamics.

In alignment with \cite{SchVee18} (and as illustrated by new results in the section ``\nameref{Sdsr}''),
it is crucial to make use of IQCs with a nontrivial terminal cost for establishing nonconservative robustness tests.
Moreover, ``\nameref{sb:iqc}'' clarifies why Theorem~\ref{Trs} with $Z\neq 0$ encompasses
the general IQC theorem in \cite{MegRan97} for rational multipliers, thus providing a tight relation between dissipativity and
IQC theory. It is remarkable that the latter requires properties on causality, stability of the uncertainty, and well-posedness of the feedback interconnection \r{fbg}, which are absent in Theorem~\ref{Trs}.

For $\Psi=I$, Corollary~\ref{Csdi} does link the existence of a positive definite solution of \r{lmif} to $A$ being Hurwitz if the left-upper block of $P$ is positive semidefinite. The latter property reads as $\Psi_1^*P\Psi_1\cge 0$ in the dynamic case. This motivates the following result,
which does not require any controllability assumptions and whose proof is identical to the one of
\cite[Th. 4]{DieSch10}.

\theorem{\label{Tiqcsta}
Let $X=X^T$ satisfy \r{lmif}, and suppose that $K\!=\!K^T$ certifies $\Psi_1^*P\Psi_1\cge 0$ as
$$
\mat{cccc}{A_\Psi&B_{\Psi_1}\\I&0}^T\!\!
\mat{cccc}{0&K&\\K&0}
\mat{cccc}{A_\Psi&B_{\Psi_1}\\I&0}\cle
\mat{ccc}{C_\Psi&D_{\Psi_1}}^T\!\!P\mat{ccc}{C_\Psi&D_{\Psi_1}}.
$$
Then, $A$ and $A_\Psi$ are Hurwitz iff
$\small\mat{cc}{X_1+K&X_{12}\\X_{12}^T&X_2}\cg 0.$
}

Since $K$ is in general indefinite, it is not true that \r{lmif}, if feasible, admits a positive definite solution (even if $A_\Psi$ and $A$ are Hurwitz).
This pinpoints a key trouble in using standard dissipativity arguments to generate robust stability tests in case that $\Psi$ involves dynamics.

\subsection{Robust performance with dynamic
integral
quadratic constraints }

Consider the uncertain interconnection
\eql{fbp}{
\mat{c}{\dot x\\z\\e}=
\mat{ccc}{
A  &B&B_d\\
C  &D&D_d\\
C_e&D_e&D_{ed}}
\mat{c}{x\\w\\d},\ \ w\in \Des(z)
}
as shown on the right in Fig.~\ref{figrp}, where $\sy$, $\sy_d$, $\sy_e$, and $\sy_{ed}$ are the sub-blocks that constitute the overall transfer matrix of the linear system in \r{fbp}. Next to the interconnection signals $(w,z)$, this loop comprises the
signals $(d,e)$ that describe a desired performance specification to hold robustly.
In general terms, $d$ is an external disturbance that excites the loop \r{fbp}, while the output $e$ is interpreted as
an error signal that should be kept small. Precisely, the target is to guarantee that there exists some $\ga_p>0$ for which the hard static IQC
\eql{rqp}{
\int_{0}^T
\mat{c}{e(t)\\d(t)}^T\!\!\!
P_p
\mat{c}{e(t)\\d(t)}\,dt=
\int_{0}^T
s_{P_p}(d(t),e(t))
\,dt\leq \ga_p\|x(0)\|^2\te{for all}T\geq0
}
holds for all trajectories of the interconnection \r{fbp}. Here, $P_p$ is a symmetric performance matrix that is partitioned into blocks
of dimensions compatible  with the signals $e$, $d$ and denoted as
$$
P_p=\mat{cc}{Q_p&S_p\\S_p^T&R_p}\te{with}Q_p\cge 0.
$$
For example, the choice $Q_p=I$, $S_p=0$, and $R_p=-\ga^2I$ captures the fact that the energy gain of the channel $d\mapsto e$ in \r{fbp}
is robustly bounded by $\ga>0$. Robust passivity of $d\mapsto e$ is characterized by taking $P_p=-P_{\pas}$.

Robust performance is certified by the dissipativity FDI
\eql{fdirp}{
\mat{cc}{\sy&\sy_d\\I&0}^*\!\!\Psi^*\mi\Psi\mat{cc}{\sy&\sy_d\\I&0}+
\mat{cc}{\sy_e&\sy_{ed}\\0&I}^*\!\!P_p\mat{cc}{\sy_e&\sy_{ed}\\0&I}
\cl0,
}
for some $P\in\bm{P}$. The series interconnection of the filter \r{psi} with the linear system in \r{fbp} leads to the matrices for expressing the corresponding dissipation LMI, which reads as
\mul{\label{lmifp}
(\pl)^T
\mat{cc|cc}{0&0&X_{1}&X_{12}\\0&0&X_{12}^T&X_{2}\hl X_{1}&X_{12}&0&0\\X_{12}^T&X_{2}&0&0}
\mat{cc|cc}{A_\Psi&B_{\Psi_1}C&B_{\Psi_1}D+B_{\Psi_2}&B_{\Psi_1}D_{d}\\0&A&B&B_d\hl I&0&0&0\\0&I&0&0}+\\
+
(\pl)^TP
\mat{cc|cc}{C_\Psi&D_{\Psi_1}C&D_{\Psi_1}D+D_{\Psi_2}&D_{\Psi_1}D_{d}}+
(\pl)^T
P_p
\mat{cc|cc}{0&C_e&D_{e}&D_{ed}\\0&0&0&I}\cl 0.
}
The placeholder $(\pl)$ is used for the corresponding matrix on the right to save space.


\theorem{\label{Trp}
Let $\Des$ satisfy a finite-horizon IQC with terminal cost for the multiplier class $\Psi^*\bm{P}\Psi$. If there exists some $P\in\bm{P}$, a corresponding terminal cost matrix $Z$, and some $X$ with \r{pos} and \r{lmifp}, then $E(X,Z):=X_2-X_{12}^T(X_1+Z)^{-1}X_{12}\cg 0$ and
all trajectories of \r{fbp} satisfy
\eql{rp}{
x(T)^TE(X,Z)x(T)+\int_0^Ts_{P_p}(d(t),e(t))\, dt
\leq x(0)^TX_{2}x(0)\te{for all}T\geq 0.
}
Therefore, robust quadratic performance \r{rqp} holds for all loop-trajectories. Moreover, the feedback loop \r{fbp} for $d=0$ [or equivalently, the loop \r{fbg}] is robustly stable.
}

\remark{\label{Rinc}
Note  that the solution set of the robust performance LMI \r{lmifp} is contained in that of \r{lmif} for robust stability. Indeed,
\r{lmif} follows  from \r{lmifp} by canceling the last block row and column and using $Q_p\cge 0$.}

\proo{As in Theorem~\ref{Trs}, the LMI \r{lmifp} guarantees that trajectories of the loop \r{fbp} filtered with \r{psi}
satisfy an SDI for the supply rate $s_{P_p}$ and a quadratic storage function defined with the matrix in \r{pos}. Then, \r{rp} is obtained by exploiting \r{pos} and noting
$$
x^T\left[X_2-X_{12}^T(X_1+Z)^{-1}X_{12}\right]x\leq \mat{c}{\xi\\x}^T\!\!\mat{cc}{X_{1}+Z&X_{12}\\X_{12}^T&X_{2}}\mat{c}{\xi\\x}
$$
for all real vectors $\xi$ and $x$ of compatible dimensions. The conclusion about robust stability
is a consequence of Remark~\ref{Rinc} and Theorem~\ref{Trs}.
}

The specialization to $\Psi=I$ is fully covered as well, which is a standard result from dissipativity theory
that can be traced to  \cite[Theorem 5]{Wil72a}. Although not difficult, only reasons of space prevent addressing more advanced performance
criteria expressed through hard IQCs on the performance channel by families of dynamic multipliers. More details can be found in \cite{VeeSch16} and references therein.

\subsection[Specialization: Input-to-state performance]{Specialization: Input-to-state performance}\label{Sinv}

It is important to not overlook specializations of Theorem~\ref{Trp} that are of interest by themselves.
For example, with $Q_p=0,\,S_p=0$, and $R_p\!=\!-\ga^2I$, the FDI \r{fdirp} simplifies to
\eql{fdirpi}{
\mat{cc}{\sy&\sy_d\\I&0}^*\!\!\Psi^*\mi\Psi\mat{cc}{\sy&\sy_d\\I&0}-\ga^2\mat{cc}{0&I}^*\mat{cc}{0&I}\cl 0.
}
This induces a related simplification of the LMI \r{lmifp}. Then, next to robust stability, Theorem~\ref{Trp} guarantees the robust input-to-state performance condition
$$x(T)^TE(X,Z)x(T)\leq
x(0)^TX_{2}x(0)+\ga^2\int_0^Td(t)^Td(t)\, dt\te{for all}T\geq 0.$$
The additional LMI constraint
\eql{lmirpi}{
\mat{ccc}{Y&0&C_e\\0&X_1+Z&X_{12}\\C_e^T&X_{12}^T&X_2}\cg 0
}
implies  $0\cle C_e Y^{-1}C_e^T\cl E(X,Z)$ (Lemma~\ref{Lschur}). Therefore, feasibility of \r{lmifp} and \r{lmirpi}
robustly guarantees, for all disturbances $d$ of finite energy, the ellipsoidal invariance property
\eql{inv}{
C_ex(T)\in \{e\in\R^{\dim(e)}\st e^TY^{-1}e\leq x(0)^TX_{2}x(0)+\ga^2\|d\|_2^2\}\te{for all}T\geq 0.
}
Minimizing the trace of $Y$ is a measure for minimizing the size of the ellipsoid on the right.
This is a special case of more consequences for such input-to-state performance guarantees if using dynamic IQCs, as addressed in \cite{Bal02,FetSch17c}.

\subsection{Specialization: Loop without external disturbance}

Similarly, the disturbance input $d$ can be absent (or set to zero), such that \r{fbg} is only excited by nonzero initial conditions. Then, the robust performance FDI \r{fdirp} simplifies to
\eql{fdirp2}{
\mat{cc}{\sy\\I}^*\!\!\Psi^*P\Psi\mat{cc}{\sy\\I}+\sy_{e}^*Q_p\sy_{e}\cl 0,
}
with a similar observation for the LMI \r{lmifp}. Next to robust stability, Theorem~\ref{Trp} then permits to robustly guarantee the state-to-output performance condition
$$x(T)^TE(X,Z)x(T)+\int_0^T e(t)^TQ_pe(t)\,dt\leq x(0)^TX_{2}x(0)\te{for all}T\geq 0.$$

\subsection[Disturbance inputs in stability analysis]{Disturbance inputs in stability analysis}\label{Siqc}

One encounters a variety of configurations and formulations concerning the characterization and
analysis of stability in the literature. So far, robust stability was confined to the specific property
\r{sta} for loops excited through nontrivial initial conditions. Note that the derived robust stability test
has wider implications. Indeed, the related LMI \r{lmif} implies the validity of the robust performance LMI \r{lmifp} for a multitude of other configurations.

This is only illustrated for the specific loop \r{fbi} related to the classical IQC theorem, as on the left in Fig.~\ref{Frs}.
Since the initial condition in \r{fbi} is taken to be zero, robust stability means to guarantee, for some $\ga>0$,
the robust performance condition
\eql{rpiqc}{
\frac{1}{\ga}\int_0^T \|z(t)\|^2\,dt\leq\ga \int_0^T\|d(t)\|^2\,dt\te{for all}T\geq 0}
on $d\mapsto z$. If matching \r{fbi} with \r{fbp} and \r{rpiqc} with \r{rqp}, the resulting FDI
\r{fdirp} reads as
\eql{fdirs}{
\mat{cc}{\sy&I\\I&0}^*\!\!\Psi^*P\Psi\mat{cc}{\sy&I\\I&0}+\mat{cc}{\sy&I\\0&I}^*\!\!\mat{cc}{\frac{1}{\ga}I&0\\0&-\ga  I}\mat{cc}{\sy&I\\0&I}\cl 0.
}

\lemma{\label{Lrs}If $X=X^T$ satisfies the robust stability LMI related to the FDI \r{fdis}, then there exists some $\ga>0$ such that the same storage function matrix $X$ satisfies the
dissipation LMI for \r{fdirs}.}
\proo{To lighten the notation, it is shown that \r{fdis} implies the validity of \r{fdirs} for some $\ga>0$ under
the assumption that $A_\Psi$, $A$ have no eigenvalues in $i\R$. The same argument applies for the related LMIs,
without any assumptions and only involving larger matrices.
Note that \r{fdirs} can be expressed  as
\eql{fdih}{
\mat{cc}{
(\Psi_1\sy+\Psi_2)^*P(\Psi_1\sy+\Psi_2)+\frac{1}{\ga}\sy^*\sy&\ \ (\Psi_1\sy+\Psi_2)^*P\Psi_1+\frac{1}{\ga}\sy^*\\
\Psi_1^*P(\Psi_1\sy+\Psi_2)+\frac{1}{\ga}\sy&\Psi_1^*P\Psi_1+\frac{1}{\ga}I-\ga I}
\cl 0.
}
All transfer matrix blocks are bounded on the compact set $i\R\cup\{\infty\}$.
By \r{fdis}, the left upper block is negative definite on $i\R\cup\{\infty\}$ for all large $\ga>0$.
The entry $-\ga I$ in the right lower block renders \r{fdih} valid if $\ga>0$ is large enough (Lemma~\ref{Lschur}).
}

Such arguments expand equally well, for example, to the characterization of robust stability for the channels
$d\mapsto\col(z,w)$ and $\col(d_1,d_2)\mapsto\col(z,u)$ in Fig.~\ref{Frs}.

\subsection[The benefit of dynamic
integral
quadratic constraints: An example]{The benefit of dynamic IQCs: An Example}
\label{Sex1}

The benefit of using dynamics in the filter $\Psi$ is illustrated by means of a detailed example.
For $\al\in[0,20]$, consider the system and controller with transfer functions
$G_\al(s)=-\frac{\al}{(s+1)(s+2)}$ and $K(s)=1$ placed in the loop on the right in Fig.~\ref{fig1}.
With a minimal realization of $G_\al$, this loop admits the description
\eql{ex1}{
\dot{x}=\mat{cc}{-3&-2\\1&0}x+\mat{c}{\al\\0}(w+d),\ \ z=e=\mat{cc}{0&-1}x,\te{and}w=\sat(z).
}

Note that the slope of the saturation $\vp=\sat$ is restricted to $[0,1]$, since
\eql{slo}{
\vp(0)=0\te{and}0\leq \frac{\vp(y)-\vp(z)}{y-z}\leq 1\te{for all}y,z\in\R,\ y\neq z.
}
This motivates embedding \r{ex1} into an interconnection with the slope-restricted uncertainty
\eql{slou}{
\Denl(z):=\{\vp(z)\st \vp:\R\to\R\text{\ satisfies \r{slo}}\} \te{for} z\in\Se.}
Recall that $\Denl(z)$ satisfies a hard static IQC for the multiplier $P_{0,1}$. Since any $\vp$ with \r{slo} is monotone, the primitive $f(z):=\int_0^z\vp(x)\,dx$ is convex with a
minimum at $0$ and its subgradient satisfies $\partial f(z)=\vp(z)$ for all $z\in\R$.
Therefore, ``\nameref{sb:zf}'' identifies a further valid dynamic IQC with terminal cost $0$.
A corresponding multiplier is given by $\Psi_{h}^*P_{\pas}\Psi_{h}$, with the filter $\Psi_h$ in
\r{zfm} determined on the basis of an impulse response function $h$ with the properties in Lemma~\ref{Lzf}.
With any $a>0$, we choose $h_a(t)=ae^{-at}\geq 0$ for $t\geq 0$ and note that $\int_0^\infty h_a(t)\,dt=1$.

This dynamic Zames-Falb IQC can be conically combined with the static one for  $P_{0,1}$ on the basis of the following instrumental observation.

\lemma{\label{Lcom}Suppose that $\Des$ satisfies an  IQC for $\Psi_a^*P_a\Psi_a$, $\Psi_b^*P_b\Psi_b$ with terminal cost matrices $Z_a$, $Z_b$, respectively. With arbitrary $\la_a,\la_b\geq 0$,  the uncertainty $\Des$ then also satisfies a hard IQC for the conically combined multiplier
$$
\la_a\Psi_a^*P_a\Psi_a+\la_b\Psi_b^*P_b\Psi_b=
\mat{c}{\Psi_a\\\Psi_b}^*\mat{cc}{\la_aP_a&0\\0&\la_bP_b}\mat{c}{\Psi_a\\\Psi_b}\te{and terminal cost}
\mat{cc}{\la_aZ_a&0\\0&\la_bZ_b}.
$$
If $\bm{P_a}$, $\bm{P_b}$ are convex cones, the combination of
$\Psi_a^*\bm{P_a}\Psi_a$, $\Psi_b^*\bm{P_b}\Psi_b$
equals
$\Psi_a^*\bm{P_a}\Psi_a+\Psi_b^*\bm{P_b}\Psi_b$.
}
This statement generalizes to a finite number of multiplier classes instead of just two.

\proo{The result follows by multiplying the IQC \r{iqcZ} for $(\Psi_a,P_a,Z_a)$, $(\Psi_b,P_b,Z_b)$ with the scalars $\la_a,\la_b\geq 0$, respectively,
and summing the two inequalities. 
}

If continuing with the example, one concludes that  $\Denl$ satisfies an IQC with terminal cost matrix $Z=0$ for the multiplier class $\Psi_{\rm zf}^*\bm{P_{\rm zf}}\Psi_{\rm zf}$ defined by
\eql{ex1mul}{
\Psi_{\rm zf}
=\mas{c?cc}{A_{\Psi_{{\rm zf}}}&B_{\Psi_{{\rm zf}}}\\\hlineB{2} C_{\Psi_{{\rm zf}}} &D_{\Psi_{{\rm zf}}}}:=
\mas{c?cccc}{-a&a&0\\\hlineB{2} 0&1&0\\0&0&1\hl -1&1&0\\ 0&0&1 }\text{\ and\ }
\bm{P_{\rm zf}}=\left\{\mat{cc}{\la_1P_{0,1}&0\\0&\la_2P_{\pas}}\st \la_1\geq 0,\ \la_2\geq 0\right\}.
}
It is shown in the section ``\nameref{Sinv}'' how to determine guaranteed bounds on
$\sup_{t\geq 0}|e(t)|$ if $x(0)=0$ and for all disturbances $d\in\Se$ with $\|d\|_2\leq 1$.
If $X$ and $Y$ satisfy \r{lmifp} and \r{lmirpi} with $Z=0$ and $C_e=\mat{cc}{0&-1}$,
such a bound is given by $\sqrt{{\rm trace}(Y)}=\sqrt{Y}$, and optimal bounds are obtained by minimizing $Y$ over
these LMI constraints.

The numerical results obtained for $G_\al$ with  $\al\in[5,50]$ and $a=10$
(using Yalmip \cite{Lof04} and Matlab's~\cite{ML20b} LMI solver) are depicted in Fig.~\ref{figex1}.
Based on $\{\la P_{0,1}\st\la\geq0\}$ (static multipliers) or $\Psi_{\rm zf}^*\bm{P_{\rm zf}}\Psi_{\rm zf}$ (dynamic multipliers), we obtain the blue and red curves, respectively.
The vertical dotted lines indicate the robust stability margins, those values of $\al$ beyond which the LMIs are infeasible
and robust stability cannot be assured. Despite their simple dynamics, Zames-Falb multipliers substantially improve the computed bounds and margins. The black curve depicts the exact amplitude bound for $e$ if replacing $\sat$ in \r{ex1} by $0$.
This is a lower bound for all robust bounds obtained with any computational technique whatsoever, since $0\in\bm{\Denl}(z)$ for all $z\in\Se$. In this example, the distance of the red and black curves can be
pushed to zero if moving $a$ to larger values, as illustrated by the green curve obtained for $a=100$.
Hence, Zames-Falb multipliers involve no conservatism.
More such experiments and comparisons with the literature can be found in \cite{FetSch17b,FetSch17c}.

Only causal Zames-Falb multipliers with a very simple parameterization are addressed in this article.
For the full class, general computational tests for standard IQCs have been developed in \cite{CheWen96,KulSaf02},
for example, and carry over to the current setting.
Multivalued nonlinearities have a wide range of applications as seen, for example, in \cite{BroLoz19,BroTan20}. It is remarkable that general Zames-Falb multipliers can be employed for subgradient nonlinearities as well. In parallel to this section,
\cite{SchHol18} contains an example that reveals their benefit over widely-used passivity-based techniques as in \cite{Bro04}.

\remark{
The classical IQC Theorem~\ref{Tiqc} does not generate these invariance results.
Still, it allows for determining the robust stability margin by finding the largest $\al\in[5,50]$
such that the FDI \r{fdis} is satisfied for some $P\in\bm{P_{\rm zf}}$.
In line with Theorem~\ref{Tiqc2}, we now show in detail why these margins are identical to the computed ones.

Let $P\in\bm{P_{\rm zf}}$ be taken such that the FDI \r{fdis} holds.
Then, \r{fdirpi} is valid for some large $\ga>0$, and the corresponding robust performance LMI has a solution $X=X^T$.
By Remark~\ref{Rinc}, $X$ also satisfies the robust stability LMI \r{lmif}.
Note that $A$ in \r{ex1} and $A_\Psi$ in \r{ex1mul} are Hurwitz, and $(C_\Psi\ \, D_{\Psi_{{\rm zf},1}})^TP(C_\Psi\ \, D_{\Psi_{{\rm zf},1}})=0$ holds. Therefore, $K=0$ certifies the dissipation LMI for $\Psi_{{\rm zf},1}^*P\Psi_{{\rm zf},1}\cge 0$ in Theorem~\ref{Tiqcsta}, which implies that $X\cg 0$. Therefore, one can take $Y$ sufficiently large to also render \r{lmirpi} for the terminal cost matrix $Z=0$ satisfied (Lemma~\ref{Lschur}).
In summary, if \r{fdis} is valid for some $P\in\bm{P_{\rm zf}}$, then the LMIs corresponding to \r{fdirpi} and \r{lmirpi} do have solutions $X$, $Y$ with $Z=0$, which proves the claim.  \epro}

\section[Computational robustness tests for structured uncertainties]{Computational robustness tests for structured uncertainties}\label{RSs}

Robust stability analysis of a complex interconnection
investigates whether it stays stable, even if some of the constituent
subsystems of such a network are affected by (possibly large)
perturbations. The results given so far are an instantiation
of this general view. However, they only involve one monolithic uncertainty and one
system in a simple feedback loop.
This section reveals that the real power of Theorems~\ref{Trs} and \ref{Trp} unfolds by allowing
for structure in the uncertainty $\Des$.

\subsection{Structured uncertainties}\label{Ssun}

Instead of an abstract description, we illustrate the emergence of structured uncertainties with an example of a simple system that is affected by several uncertain components of a different nature in a structured fashion.

With uncertain parameters $\delta_1\in\Se$, $\delta_2\in\R$ and the nonlinear element $\sat$, consider
\eql{ex2a}{
\dot x_1=-x_1+\frac{1}{1+\delta_1}x_2,\ \ \dot x_2=\delta_1x_1-x_2+\delta_2\sat(x_2-2x_1).
}
By introducing auxiliary signals, this system can be equivalently expressed as
\eql{ex2lfr}{
\arraycolsep.5ex
\left.\arr{rcl}{
\mat{c}{\dot x_1\\\dot x_2}&=&
\mat{rr}{-1&1\\0&-1}
\mat{c}{x_1\\x_2}+
\mat{rr|r|r}{-1&0&0&0\\0&1&1&1}
\mat{c}{w_1\\ w_2\hl w_3\hl  w_4}
\\
\arraycolsep.5ex
\mat{c}{z_1\\ z_2\hl z_3\hl z_4}&=&
\mat{rr}{0&1\\1&0\hl 0&0\hl -2&1}x+
\mat{rr|r|r}{-1&0&0&0\\0&0&0&0\hl 0&0&0&1\hl 0&0&0&0}
\mat{c}{w_1\\ w_2\hl w_3\hl  w_4}
}\!\right\}\!,\
\mat{c}{w_1\\ w_2\hl w_3\hl w_4}=
\mat{cc}{\delta_1z_1\\\delta_1z_2\hl \delta_2z_3\hl \sat(z_4)}.
}
Note that \r{ex2lfr} is a feedback loop of a linear system and a static nonlinearity with a diagonal structure.
In robust control, this is called a linear fractional representation, and Matlab's robust control toolbox generates such descriptions with ease.
More details about this perspective can be found in \cite{ZhoDoy96}.
Note that the generation of such representations is in full alignment with
Willems' approach of tearing, zooming, and linking to modeling in the behavioral framework \cite{Wil07csm}.
Despite its now rather long history and great success, the power of
such structured representations still seems not fully appreciated in the control community as a whole.

Now suppose that the parameters are bounded as $|\delta_1(t)|\leq r$,  for $t\geq 0$ and $|\delta_2|\leq 1$, and
recall that $\sat$ is captured by \r{slou}. This motivates the introduction of the individual uncertainties
\gan{
\bm{\Delta_a}(z_a):=\{\delta z_a\st \delta\in\Se,\ |\delta(t)|\leq r\text{\ for\ }t\geq 0\}\te{and}
\bm{\Delta_b}(z_b):=\{\delta z_b\st \delta\in\R,\ |\delta|\leq 1\}}
for $z_a\in\Se^2$, $z_b\in\Se$ and $\bm{\Delta_c}:=\Denl$. Then, \r{ex2a} is viewed as an uncertain interconnection with the
overall uncertainty
$
\Des(z)=\left\{\col(w_a,w_b,w_c)\st w_a\in\bm{\Del_a}(z_a),\ w_b\in\bm{\Del_b}(z_b),\ w_c\in\bm{\Del_c}(z_c)\right\}
$. 

The individual uncertainties satisfy IQCs with terminal cost $0$ for, respectively,
\gan{
\bm{P_a}:=\left\{\mat{cc}{r^2Q&rS\\rS^T&-Q}\st Q,S\in\R^{2\times 2},\ Q\cge 0,\ S^T\!+\!S=0\right\},\
\bm{P_b}:=\left\{\la\mat{cc}{1&0\\0&-1}\st \la\geq 0\right\}
}
and \r{ex1mul}. Indeed, the IQC for $\bm{\Del_a}$ and $P\in\bm{P_a}$ follows, since $|\delta(t)|\leq r$ implies
$s_P(\delta(t) I_2,I_2)=r^2Q-Q|\delta(t)|^2+r(S^T+S)\delta(t)=
Q(r^2-|\delta_2(t)|^2)\cge 0$ for all $t\geq 0$. For $\bm{\Del_b}$, note that $z\mapsto \delta z$ with $|\delta|\leq 1$ also satisfies a sector condition for the multiplier $P_{-1,1}$.

As in Lemma~\ref{Lcom} and reminiscent of the classical S-procedure \cite{GusLik06},
it is straightforward to diagonally combine these individual IQCs to one for the overall uncertainty
$\Des$. We state the general fact for two uncertainties to simplify notations.

\lemma{\label{Lcom2}Suppose that $\bm{\Des_a}$ and $\bm{\Des_b}$ satisfy IQCs for the multipliers $\Psi_a^*P_a\Psi_a$ and $\Psi_b^*P_b\Psi_b$ with symmetric terminal cost matrices $Z_a$ and $Z_b$, respectively. Then, the diagonal combination
$$
\Des(z_a,z_b):=\left\{\col(w_a,w_b)\st w_a\in\bm{\Del_a}(z_a),\ w_b\in\bm{\Del_b}(z_b)\right\}
$$
satisfies an IQC for the multiplier $\Psi^*P\Psi$ with the terminal cost matrix $Z$, where
$$\arraycolsep.3ex
\Psi=\mat{cccc}{\Psi_{a,1}&0&\Psi_{a,2}&0\\0&\Psi_{b,1}&0&\Psi_{b,2}},\ P=\mat{cc}{\la_aP_a&0\\0&\la_bP_b},\text{ and }
Z=\mat{cc}{\la_aZ_a&0\\0&\la_bZ_b}\text{\ with\ }\la_a,\la_b\geq 0.
$$
The column partitions $\Psi_a=\mat{cc}{\Psi_{a,1}&\Psi_{a,2}}$, $\Psi_b=\mat{cc}{\Psi_{b,1}&\Psi_{b,2}}$ are taken according to the row partitions of $\col(z_a,w_a)$, $\col(z_b,w_b)$, respectively.
If $\bm{P_a}$, $\bm{P_b}$ are convex cones, the combination of
$\Psi_a^*\bm{P_a}\Psi_a$ and $\Psi_b^*\bm{P_a}\Psi_b$ leads to
$\Psi^*\bm{P}\Psi$ with $\bm{P}=\{\diag(P_a,P_b)\st P_a\in\bm{P_a},\ P_b\in\bm{P_b}\}$.
}

In continuing the example, $\bm{P}$ consists of all matrices $\diag(P_a,P_b,P_c)$ with $P_a\in\bm{P_a}$, $P_b\in\bm{P_b}$,
and $P_c\in\bm{P_c}$, while $\Psi$ is composed of
$$\arraycolsep.4ex
\mat{c|c}{\Psi_{a,1}&\Psi_{a,2}}=\mat{c|c}{1&0\\0&1},\
\mat{c|c}{\Psi_{b,1}&\Psi_{b,2}}=\mat{cc|cc}{1&0&0&0\\0&1&0&0\\0&0&1&0\\0&0&0&1},\
\mat{c|c}{\Psi_{c,1}(s)&\Psi_{c,2}(s)}=\mat{c|c}{ 1&0\\0&1 \\ 1-\frac{a}{s+a}&0\\0&1 }.
$$

Despite its simplicity, this example demonstrates all the ingredients for generating tailored robustness test
for complex interconnections that are computational in general. It shows how to build representations \r{fbg} of uncertain systems
with a structured $\Des$ that collects multiple uncertainties of a different nature in a systematic fashion. Moreover, it illustrates
how to exploit information about individual uncertainties to choose suitable individual multiplier classes, and how they are
combined into a multiplier class for the overall uncertainty.

This strategy is modular. One can strengthen the test by choosing more powerful (larger) multiplier classes for the individual uncertainties, at the expense of a higher computational cost. In our example, one might use  better Zames-Falb multipliers for the
static nonlinearity to reduce conservatism. Similarly, one can reduce the computational cost by coarsening the set of individual multipliers, if accepting a potential  increase in conservatism.

The construction of relaxation hierarchies of such families is addressed, for example, with matrix sum-of-squares techniques in \cite{Sch05,SchHol06,VeeSch16}. In \cite{VeeSch16} and references therein, it is shown how to employ structured dynamic multipliers for performance specifications, which is not detailed any further in this article.
Despite all the theoretical progress in past years, the development of customized LMI solvers
for increasing computational efficiency, for example by exploiting the described modularity,
is highly desirable but still in its infancy \cite{WalKao08}.

\subsection{Nonlinear uncertainties}

Dissipativity arguments generate IQCs for dynamic uncertainties.
For example, consider \r{nlsy} and suppose that this system has the equilibrium $(0,0)$ and is $s_P$-dissipative with a nonnegative
storage function $V$ satisfying $V(0)=0$. The dissipation inequality implies that
\eql{dnl}{\Dednl(z):=\{w\in\Se^\ny\st x,\dot x\in\Se^\nx,\ \dot x=f(x,z),\ w=g(x,z),\ x(0)=0\}\te{for}z\in\Se^\nuu}
satisfies a hard IQC for $\bm{P}:=\{\la P\st\la\geq 0\}$. If considering \r{ex2a} with this nonlinearity
instead of $\sat$, one generates a corresponding robustness test by replacing
$(\Denl,\Psi_{\rm zf}^*\bm{P_{\rm zf}}\Psi_{\rm zf})$ with $(\Dednl,\bm{P})$.
For robust loop stability, one should note that strict dissipativity
or additional assumptions on the storage function $V$ permit to draw conclusions about the ``internal'' state-trajectory
in $\Dednl(z)$ as well, which was emphasized in \cite{Wil72a} and is addressed in detail in \cite{ArcMei16}.

The compositional framework in \cite{ArcMei16} with many inspiring examples builds on \cite{Wil72a,MoyHil78}. It is subsumed to \r{fbp} with a nondynamic linear system and a diagonally structured $\Des$, as in Lemma~\ref{Lcom2},
possibly involving a large number of blocks instead of just two. The linear system then defines a static interconnection of the dynamical components in $\Des$.
In our terminology, the main result of \cite{ArcMei16} relies on a class of dynamic multipliers $\Psi^*\bm{P}\Psi$,
composed as in Lemma~\ref{Lcom2}, for which $\Des$ satisfies a dynamic IQC with terminal cost $Z=0$. Then, \cite[Proposition 8.2]{ArcMei16} guarantees interconnection dissipativity through the existence of a certificate $X\cge 0$ of the FDI
$$
\mat{cc}{D&D_d\\I&0}^*\!\!\Psi^*\mi\Psi\mat{cc}{D&D_d\\I&0}+
\mat{cc}{D_e&D_{ed}\\0&I}^*\!\!P_p\mat{cc}{D_e&D_{ed}\\0&I}\cle0\te{for some}P\in\bm{P}.
$$
Working with Theorem~\ref{Trp} involves the use of nontrivial terminal cost matrices $Z$ and the coupled positivity constraint $X+Z\cge 0$. This not only generalizes the main result of \cite{ArcMei16}, but also
clarifies the reason for the discrepancies to classical IQC theory, as elaborated on in \cite[Section 9.3]{ArcMei16}.

If, next to the system, the uncertainty itself also comprises only static sector-bounded maps, \r{fbp} defines a static
Lur'e system, which covers deep neural networks \cite{SuyVan99},  for example. Theorem~\ref{Trp} in combination with
the geometric interpretation of classical Lagrange relaxations for indefinite quadratic functions \cite{BoyGha94,BoyVan04},
then permits to construct computational safety guarantees as (for example) recently addressed in  \cite{FazMor19}. The range of further applications is stunningly wide.

\subsection{Dynamic uncertainties and static supply rates}

Capturing unmodeled dynamics of possibly high order is of key relevance in robust control. For example, a damping coefficient in a mechanical system might result from ignoring dynamical effects in a passive shock absorber mechanism.

Concretely, consider an uncertainty defined with stable transfer functions $\delta$ of any degree.
These are viewed as single-input, single-output systems described by $\hat w(s)=\delta(s)\hat z(s)$ in the frequency domain, where $\hat{\phantom{a}}$ denotes the Laplace transform. With a minimal realization $\delta(s)=C_\delta(sI-A_\delta)^{-1}B_\delta+D_\delta$,
the matrix $A_\delta$ is Hurwitz and the related system is described by
\eql{syu}{\dot x=A_\delta x+B_\delta z,\ \ w=C_\delta x+D_\delta z,\ \ x(0)=0.}
This is shortly expressed as $w=\delta z$. Now suppose that the transfer function $\delta$ is positive real, which translates into the FDI $s_{P_{\pas}}(\delta,1)=\delta^*+\delta\geq 0$. By Corollary~\ref{Cdi} and since $A_\delta$ is Hurwitz,  there exists some $\Xs\cge 0$ with $x(T)^T\Xs x(T)\leq \int_0^T s_{P_{\pas}}(w(t),z(t))\,dt$ for all $T\geq 0$ and all trajectories of \r{syu}. Hence, $\delta$ satisfies a hard IQC for the multiplier $P_{\pas}$.

This example demonstrates the general principle of how dissipativity generates valid IQC for dynamic uncertainties of practical interest. Again, more information about $\delta$ increases the set of multipliers and typically reduces conservatism of the related robustness test. For example, it may be known that the $H_\infty$-norm of $\delta$ is bounded as $\|\delta\|_\infty:=\sup_{\om\in\R}|\delta(\io)|\leq 1$.
Again by Corollary~\ref{Cdi}, the corresponding FDI $s_{P_{\gai}}(\delta,1)=1-\delta^*\delta\geq 0$ translates into a hard IQC for $P_{\gai}$. By Lemma~\ref{Lcom},  $\delta$ also satisfies a hard IQC for the multiplier class $\bm{P}:=\{\la_1P_{\pas}+\la_2P_{\gai}\st\la_1\geq 0,\ \la_2\geq0\}$. The resulting robust stability or performance tests are superior to separate small-gain or passivity criteria \cite{GriAnd07}. However, this is not often exploited for dynamic uncertainties in practice.

As addressed for  \r{ex2a}, $\bm{P}$ is also a valid multiplier class for time-invariant and time-varying
parametric uncertainties that satisfy $\delta\in[0,1]$ or $\delta(t)\in[0,1]$ for all $t\geq0$, respectively. This reveals the limitations of passivity- and small-gain based robustness tests, even if combined.
The question arises how to exploit the specific nature of $\delta$ to generate more powerful IQCs.

\subsection[Dynamic supply rates for dynamic uncertainties]{Dynamic supply rates for dynamic uncertainties}\label{Sdsr}

This section reveals a systematic construction of dynamic IQCs with a nonzero terminal cost for dynamic uncertainties that has not appeared in the literature. Suppose that the stable transfer function $\delta$ satisfies
\eql{fdidel}{
\mat{c}{1\\\delta}^*\!\!P_0\mat{c}{1\\\delta}\cg 0\te{with some fixed}P_0=\mat{cc}{q&s\\s&r}\in\S^2,\ r\leq 0.
}
It is technically favorable to work with strict FDIs, which means that the Nyquist plot of $\delta$ is contained in the interior of the disk or half-plane corresponding to $P_0$,
instead of in its closure. This causes no limitations for robustness analysis from a practical perspective.

Then choose a stable transfer matrix $\psi$ of dimension $l\times 1$ that has full column rank on the extended imaginary axis $i\R\cup\{\infty\}$, and take any $M\in\S^l$ such that $\psi^*M\psi\cg 0$.
If combined with \r{fdidel}, this implies with
$\Psi:=\diag(\psi,\psi)$ and using elementary Kronecker algebra that
\eqn{\arraycolsep.2ex
\mat{c}{1\\\delta }^*\!\!\Psi^*(P_0\otimes M)\Psi\!\mat{c}{1\\\delta}=
\left(\mat{c}{1\\\delta}^*\!\!\!\otimes\!\psi^*\right)\!(P_0\otimes M)\!\left(\mat{c}{1\\\delta}\!\otimes\!\psi\right)=
\mat{c}{1\\\delta}^*\!\!P_0\!\mat{c}{1\\\delta}\psi^*M\psi\cg 0.
}
Therefore, by Corollary~\ref{Cdi}, the system related to $\Psi\,\col(1,\delta)$ with any realization satisfies a strict dissipation inequality for the supply rate matrix $P_0\otimes M$. However,
the existence of a positive definite storage function cannot be assured. Hence, it cannot be guaranteed that
$\delta$ satisfies a hard IQC for the multiplier $\Psi^*(P_0\otimes M)\Psi$. It is a new result that $\delta$ does satisfy an IQC with a terminal cost that is determined by an LMI certificate for the FDI $\psi^*M\psi\cg 0$.

\theorem{\label{Tdy}Let $K=K^T$ be a certificate of $\psi^*M\psi\cg 0$. Then, the dynamic uncertainty
$$
\Dedy(z):=\{\delta z\st \delta\text{\ is a stable transfer function with \r{fdidel}}\} \te{for} z\in\Se
$$
satisfies an IQC with multiplier $\diag(\psi,\psi)^*(P_0\otimes M)\diag(\psi,\psi)$ and terminal cost $Z=P_0\otimes K$.}

The proof is provided in ``\nameref{sb:Tdy}.'' This result exhibits the essential point in Definition~\ref{DiqcZ},
to be viewed as a partial dissipativity condition that does not involve the full state in a realization of $\Psi\,\col(1,\delta)$.
Initial steps to exploit this idea for IQC robustness results involving infinite-dimensional uncertainties have been taken in \cite{BarSch20}.

Theorem~\ref{Trp} leads to the following robust performance test for \r{fbp} with $\Des:=\Dedy$. Choose a minimal realization $(A_\psi,B_\psi,C_\psi,D_\psi)$ for $\psi$ (such that $A_\psi$ is Hurwitz) and the resulting one for
$\Psi=\diag(\psi,\psi)$ with $\pl_\Psi=\diag(\pl_\psi,\pl_\psi)$ for $\pl\in\{A,B,C,D\}$ (such that $A_\Psi$ is Hurwitz).

\corollary{\label{Cdy}All trajectories of the loop \r{fbp} with $\Des:=\Dedy$ satisfy \r{rp} if there exist $X=X^T$, $K=K^T$ and $M\in\S^l$
such that $P=P_0\otimes M$ and $Z:=P_0\otimes K$ satisfy \r{pos}, \r{lmifp} and
\eql{lmiK}{
\mat{cc}{A_\psi&B_\psi\\I&0}^T\!\!\mat{cc}{0&K\\K&0}\mat{cc}{A_\psi&B_\psi\\I&0}
\cl \mat{cc}{C_\psi&D_\psi}^T\!\!M\mat{cc}{C_\psi&D_\psi}.
}}

Corollary~\ref{Cdy} involves a set of LMI constraints on $X$, $K$, $M$, and, hence, also on the terminal cost matrix $Z$. This allows for convex optimization over all these variables, for example, if including further dissipativity constraints
as in the section ``\nameref{Sinv}.'' This constitutes
a substantial progress over the generic choice of $Z$ in Theorem~\ref{Tiqc2}, which involves a nonlinear coupling of $Z$ and $P$
that cannot be convexified. In contrast to the terminal cost matrix proposed in \cite{FetSch17c} for general IQCs, the newly suggested one is tight. This follows if the LMIs in Corollary~\ref{Cdy} are feasible, since
$K$ can be taken such that $P_0\otimes K$ is arbitrarily close to the stabilizing solution of the related algebraic Riccati equation \r{are}.
Note that \cite{Bal02} proposes to work with $M\cg 0$ and $K=0$,
which is shown to be tight only under the restrictive assumption that $(C_\psi\ \,D_\psi)$ has full column rank.

\example{\label{Ers}
For a numerical example, consider the variant $\dot x_1=-x_1+x_2$, $\dot x_2=-x_2+\al\de(x_2-2x_1)$
of \r{ex2a} for $\delta_1=0$, a fixed parameter $\delta_2=\al\in[0.9,1]$ and a dynamic uncertainty $\delta$ replacing $\sat$. This is represented as $z=G_\al w$, $w=\de z$, with $G_\al$ described by
$$
\dot x=\mat{rrr}{-1&1 \\0&-1}x+\mat{c}{0\\\al}w,\ \ z=\mat{cc}{-2& 1}x.
$$
Place $\sy:=G_\al$ and $\Des:=\Dedy$ with $P_0:=P_{\gai}$ into the loop in Fig.~\ref{fig1} for $K=1$.
Using the performance output $G_\al e$ instead of $e$, we perform exactly the same experiment as in the section ``\nameref{Sex1},'' but now based on Corollary~\ref{Cdy}.
The guaranteed amplitude bounds are depicted in Fig.~\ref{figex3}, with the blue and red curves obtained
for the choices $\psi(s)=1$ (static multiplies) and $\psi(s)=\col(1,(s+0.9)^2/(s+1)^2)$ (dynamic multipliers), respectively.
The black curve results if taking $M\cg 0$ and $K=0$, which indicates the conservatism of \cite{Bal02} for general $\psi$'s.
Note that these results have been obtained with a fixed $\psi$ and a free parameter $M\in\S^l$ in $\psi^*M\psi$.
The choice of general parameterizations with guaranteed approximation properties is by now well understood
\cite{ChoTit99,SchKoe12} and can be exploited in the current context as well.\epro
}

Both the reduction of conservatism (as seen in the example) and the flexibility of the construction that leads to Theorem~\ref{Tdy} are viewed as key advantages. If we choose $\psi$ of dimension $l\times k$ with $k$ columns, this leads to IQCs and robustness tests for so-called repeated dynamic blocks $\delta I_\nz$. Stable  full uncertain transfer matrices $\Delta$ of dimension $m\times k$ are handled similarly, while parametric uncertainties have been addressed
in \cite[Theorem 13]{SchVee18}. By modularity, this captures all the standard types of uncertainties as covered in classical structured singular value theory, with unprecedented dissipativity proofs for the $D$ and $D/G$ scalings upper bounds  \cite{PacDoy93}. In view of Theorem~\ref{Tiqc2} and in combination with
the results in \cite{SchKoe13}, it is remarkable that these robust dissipativity tests are lossless for the so-called small block structures in $\mu$-theory \cite[Section 9]{PacDoy93}, if only the basis transfer function $\psi$ is chosen sufficiently rich.

This is a very pleasing resolution of the challenge mentioned in the introduction that was encountered
approximately 25 years ago, when working on the first version of our lectures notes about LMIs in control \cite{SchWei99}.

\section{Conclusion}

Based on the concept of dissipativity as developed by Jan Willems 50 years ago, this article surveys in a tutorial style the key principles underlying the systematic construction of robust stability and performance tests for uncertain feedback interconnections. Emphasis was placed on computational procedures and the modular composition of the corresponding LMI constraints for their implementation.

Building on the seminal work of Popov, Yakubovich, and Zames, an encompassing framework for robust stability analysis based on  IQCs has been introduced 25 years ago that relies on capturing the properties of uncertainties
through infinite horizon energy inequalities on their input-output trajectories.

The intermediate notion of finite-horizon IQCs with a terminal cost was shown to provide a seamless link between dissipativity and IQCs. New results on the systematic construction of nontrivial
classes of such IQCs for dynamic uncertainties pave the way for further extensions. Design methods
for synthesizing robust and gain-scheduled controllers based on this concept are only beginning to emerge.




\processdelayedfloats 

\newpage


\begin{figure}
\center
\scalebox{1}{\fbs{$\vp$}
\hspace{10ex}
\begin{tikzpicture}[xscale=1,yscale=1]
\node[sy3] (g) at (0,0)  {$G$};
\node[sum, left = \dl of g] (su) {\scalebox{0.7}{$+$}};
\node[sy3, left = \dl of su] (sat) {sat};
\node[sy3, left = \dl of sat] (k) {$K$};

\node[jun,right=\dl of g] (j1) {};
\node[coordinate,below = .8*\dl of g] (j2) {};
\node[coordinate,right=\dl of j1] (j3) {};

\draw[->] (k) -- node[pos=.4]{} (sat);
\draw[-] (g) -- (j1);
\draw[->] (j1)  -- node[]{$e$} (j3);
\draw[-] (j1) |- (j2);
\draw[<-] (k.west)-- ($(k.west)+(-\dl,0)$) |- (j2);
\draw[->] ($(su)+(0,1.5*\dl)$) -- node[pos=.3]{$d$} (su);
\draw[->] (sat)  -- (su);
\draw[->] (su)  -- (g);
\end{tikzpicture}}
\caption{Standard feedback loops.  The configuration for the Lur'e problem is shown on the left, while the right figure depicts a
classical saturated control loop.\label{fig1}}
\end{figure}
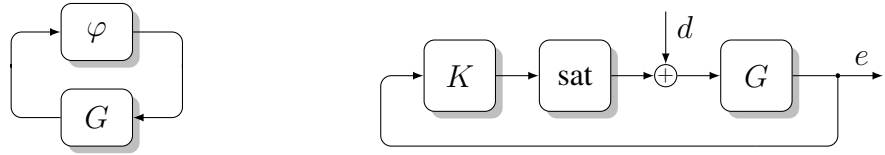
\begin{figure}
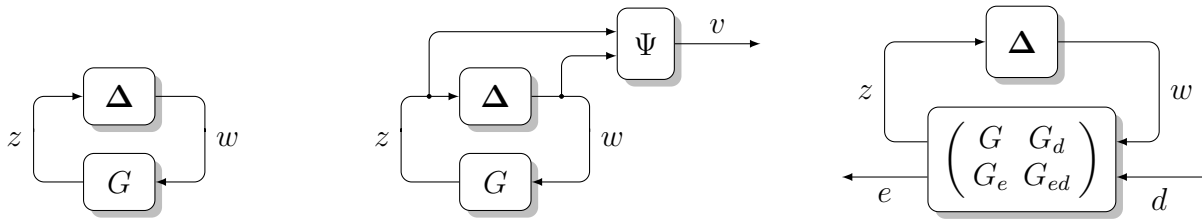

\center
\scalebox{1}{\fbp{$\Des$}}
\caption{Uncertain interconnections. These configurations are used to define robust stability (left),
analyze robust stability with integral quadratic constraints (middle), and define robust performance (right). \label{figrp}}
\end{figure}
\begin{figure}[h]
\cen{
\scalebox{1}{
\fstaba{$\Des$}\hspace{6ex}\fstabb{$\Des$}\hspace{6ex}\fstabc{$\Des$}
}
}
\caption{
Several configurations for robust stability analysis with external
disturbance inputs. As seen in the section ``\nameref{Siqc},''
such loops can be addressed in a unified fashion, on the basis of
the two main robustness results (Theorem~\ref{Trs} and Theorem~\ref{Trp}).
\label{Frs}}
\end{figure}

\mbox{}\newpage
\begin{figure}
\center
{\includegraphics[width=15cm,,trim=80 6 60 17,clip=true]{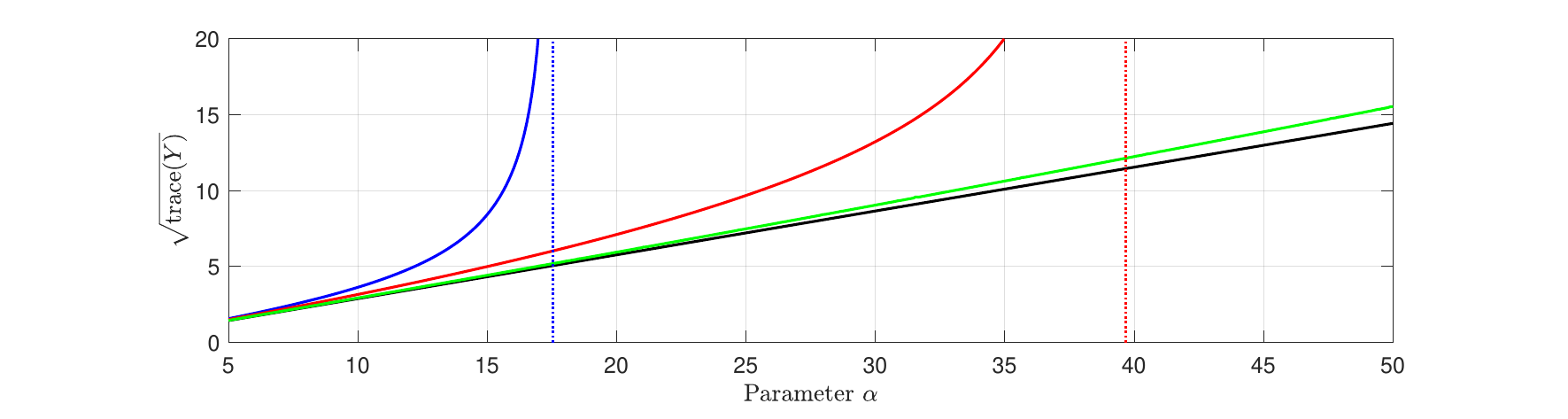}}
\caption{
Numerical results for the example in the section ``\nameref{Sex1}'' with a nonlinear uncertainty.
The graphs depict the computed guaranteed bounds on the output amplitude $\sup_{t\geq0}|e(t)|$, for the parameter values $\al\in[0,50]$ and the nominal loop (black), the uncertain loop analyzed with static multipliers (blue), and the uncertain loop handled with Zames-Falb multipliers (red for the multiplier pole $a=10$, green for the multiplier pole $a=100$). \label{figex1}}
\end{figure}

\begin{figure}[h]
\center
{\includegraphics[width=15cm,,trim=80 6 60 17,clip=true]{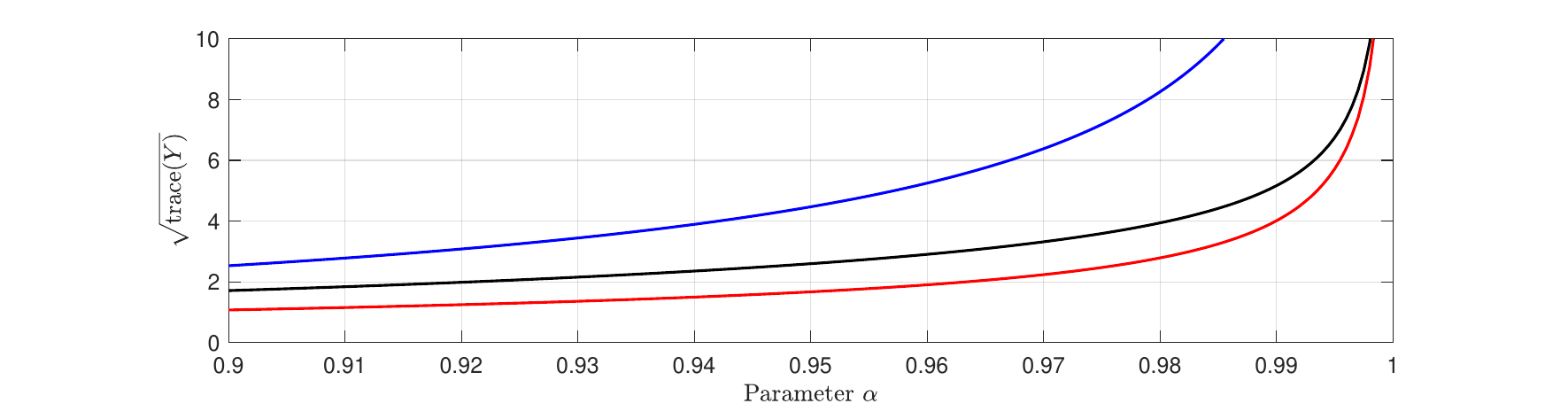}}
\caption{
Numerical results for Example \ref{Ers} with a dynamic uncertainty.
The graphs show the computed guaranteed bounds on the output amplitude $\sup_{t\geq0}|e(t)|$ for the parameter values $\al\in[0.9,1]$ and
with static multipliers (blue) and dynamic multipliers with the novel
technique  (red). A comparison with the results obtained
under the constraints in \cite{Bal02} (black) \label{figex3} reveals the conservatism of existing methods.}
\end{figure}

\sidebars 

\clearpage

\section[Article summary]
{Article summary}
\label{sb:sum}

A central notion in systems theory is dissipativity, which was introduced by Jan Willems
with the explicit goal of arriving at a fundamental understanding of the stability properties of feedback interconnections.
In robust control, the framework of integral quadratic constraints (IQCs) builds on the seminal contributions of Yakubovich and Zames in the 1960s. It provides a technique for analyzing the stability of an interconnection of some linear system in feedback with a whole class of systems, also refereed to as uncertainty.

This article surveys the key ideas of exploiting dissipativity and integral quadratic constraints to systematically construct computational tests for robust stability and performance of uncertain interconnections in terms of linear matrix inequalities. The article focuses on the recently introduced notion of finite-horizon IQCs with a terminal cost, which is shown to provide a seamless link between dissipativity theory and absolute stability theory based on dynamic IQCs.

\newpage
\section[Summary: Dissipativity for linear systems and quadratic supply rates]
{Summary: Dissipativity for linear systems and quadratic supply rates}
\label{sb:dis}

\corollary{\label{Cdi}
If the linear system \r{sy} with the transfer matrix \r{tf} is controllable and $P\in\S^{k+m}$, the following statements are equivalent:
\enu{
\item {\bf Time-domain characterization.} There exists a storage function $V:\R^n\to\R$ such that, for all trajectories of the linear system, the following dissipation inequality holds:
\eql{cdi}{
V(x(t_2))+\int_{t_1}^{t_2} \mat{c}{y(t)\\u(t)}^T\!\!P\mat{c}{y(t)\\u(t)}\,dt\leq V(x(t_1))\te{for all}0\leq t_1\leq t_2.
}
\item {\bf Frequency-domain characterization.} The following FDI holds:
$$
\mat{c}{0\\y_0}=\mat{cc}{A-\io I&B\\C&D}\mat{c}{x_0\\u_0},\  \om\in\R\ \ \Longrightarrow\ \
\mat{c}{y_0\\u_0}^*\!\!P\mat{c}{y_0\\u_0}\leq 0.
$$
This is equivalent to the FDI
$$
\mat{c}{\sy(i\om)\\I}^*\!\!P\mat{c}{\sy(i\om)\\I}\cle 0
\te{for all $\om\in \R$ such that $\io$ is no pole of $\sy$.}
$$
\item {\bf LMI characterization.} There exists a symmetric solution $\Xs$ of the dissipation LMI
\eql{cdilmi}{
\mat{cc}{A&B\\I&0}^T\!\!\mat{cc}{0&\Xs\\\Xs&0}\mat{cc}{A&B\\I&0}+
\mat{cc}{C&D\\0&I}^T\!\!P\mat{cc}{C&D\\0&I}\cle 0.
}
}
Moreover, the quadratic storage function $V(x)=x^{T}\Xs x$ with $\Xs\in\S^\nx$ satisfies the dissipation inequality \r{cdi} iff
$\Xs$ satisfies the LMI \r{cdilmi}.
Finally, if the left-upper block of $P$ is positive semidefinite and $A$ is Hurwitz, then
all solutions of \r{cdilmi} are positive semidefinite.
}

\corollary{\label{Csdi}
For a general linear system \r{sy} with the transfer matrix \r{tf} and $P\in\S^{k+m}$, the following statements are equivalent:
\enu{
\item {\bf Time-domain characterization.} There exists a storage function $V:\R^n\to\R$
and $\ve>0$ such that any system trajectory satisfies the following dissipation inequality
for $0\leq t_1\leq t_2$:
\eql{csdi}{
V(x(t_2))+\int_{t_1}^{t_2}
\mat{c}{y(t)\\u(t)}^T\!\!P\mat{c}{y(t)\\u(t)}
\,dt\leq V(x(t_1))
-\ve \int_{t_1}^{t_2}\left\|\mat{c}{x(t)\\u(t)}\right\|^2dt.
}
\item {\bf Frequency-domain characterization.} The following FDI holds:
$$\arraycolsep.2ex
\mat{c}{0\\y_0}=\mat{cc}{A-\io I&B\\C&D}\underbrace{\mat{c}{x_0\\u_0}}_{\neq 0},\ \om\in\R\te{or} y_0=Du_0,\ u_0\neq 0\ \Longrightarrow\ \
\mat{c}{y_0\\u_0}^\ast\!\!P\mat{c}{y_0\\u_0}<0.
$$
If $\eig(A)\cap i\R=\emptyset$, this is equivalent to the strict FDI
$$
\mat{c}{\sy(i\om)\\I}^*\!\!P\mat{c}{\sy(i\om)\\I}\cl 0
\te{for all}\om\in \R\cup\{\infty\}.
$$
\item {\bf LMI characterization.} There exists a symmetric solution $\Xs$ of the strict dissipation LMI
\eql{csdilmi}{
\mat{cc}{A&B\\I&0}^T\!\!\mat{cc}{0&\Xs\\\Xs&0}\mat{cc}{A&B\\I&0}+
\mat{cc}{C&D\\0&I}^T\!\!P\mat{cc}{C&D\\0&I}\cl0.
}
}
Moreover, the quadratic storage function $V(x)=x^{T}\Xs x$ with $\Xs\in\S^\nx$ satisfies the dissipation inequality \r{csdi} for some $\ve>0$  iff $\Xs$ satisfies the LMI \r{csdilmi}.
Finally, if the left-upper block of $P$ is positive semidefinite, then all solutions of \r{csdilmi} are positive definite iff $A$ is Hurwitz.
}

All throughout, the nonstrict or strict FDIs in b) are abbreviated, respectively, as
$$
\mat{c}{\sy\\I}^*\!\!P\mat{c}{\sy\\I}\cle 0\te{or}\mat{c}{\sy\\I}^*\!\!P\mat{c}{\sy\\I}\cl 0.
$$
To simplify the exposition and without harm, this convention is maintained even if $A$ does have eigenvalues in $i\R$.
Moreover, for a given realization of $\sy$, we say that $\Xs$ certifies these FDIs
(or is a certificate thereof)  if $\Xs$
satisfies the corresponding LMIs \r{cdilmi} or \r{csdilmi}, respectively.

Concrete instantiations of these dissipativity results for particular index matrices $P$
lead to strict and nonstrict versions of a variety of specializations that are often treated independently in the literature. For example, if recalling \r{P}, several formulations of the bounded real lemma are obtained with $P=P_{\gai}$, while the positive real lemma is covered with $P=-P_{\pas}$.

\newpage
\section[The circle criterion: Dissipativity proof and discussion]
{The circle criterion: Proof and Discussion}
\label{sb:pcir}

The idea of the proof of Theorem~\ref{Tcir} is as follows. Introduce the supply rates
$$
s_1(z,w)=\mat{c}{z\\w}^T\!\! P\mat{c}{z\\w}\te{and}s_2(u,y)=-\mat{c}{y\\u}^T\!\! P\mat{c}{y\\u}.
$$
Then, the inequality \r{qccir} implies that the nonlinear system $w=\nl(z)$ is dissipative with respect
to the supply rate $s_1$ (Theorem~\ref{Tddis}). Moreover, the FDI \r{fdicir} and $\eig(A)\cap i\R=\emptyset$ guarantee that the linear system \r{sy} in the loop is strictly $s_2$-dissipative with a quadratic storage function $V(x)=x^T\Xs x$ (Corollary~\ref{Csdi}).
In Willems' terminology \cite{Wil72a}, $u=w$ and $z=y$ constitutes a neutral interconnection for these two dissipative systems,
since $s_1(z,w)+s_2(u,y)=s_1(y,u)+s_2(u,y)=0$ holds for all vectors satisfying these equations. Therefore, summing the two dissipation inequalities along trajectories of the interconnection gives
\eql{h0}{
x(T)^T\Xs x(T)+\ve\int_{0}^T \|x(t)\|^2+\|u(t)\|^2\,dt \leq x(0)^T\Xs x(0)\te{for all}T\geq 0.
}
This is a simple version of one of the key results in \cite{Wil72a}, which roughly states that the
neutral interconnection of dissipative systems stays dissipative. Note that the definition of dissipativity in the current article does not impose any sign constraint on $\Xs$. The essence of the  nominal stability hypothesis in Theorem~\ref{Tcir} is to guarantee that  $\Xs$ is positive definite. Then, \r{h0} implies
\eql{h0s}{
\int_{0}^T \|x(t)\|^2+\|w(t)\|^2\,dt \leq \frac{1}{\ve}\la_{\rm max}(\Xs)\|x(0)\|^2\te{for all}T\geq 0
}
(with $\la_{\rm max}(\Xs)>0$ denoting the largest eigenvalue of $\Xs$) and, therefore, stability of the loop.


\subsection{Formal proof of Theorem~\ref{Tcir}}

By Corollary~\ref{Csdi}, the FDI \r{fdicir} and $\eig(A)\cap i\R=\emptyset$ imply that \r{csdilmi} has a solution $\Xs$. To show $X\cg 0$, apply Theorem~\ref{Tlmipos} for $-P$. Since
\r{iqcnom} implies $s_{-P}(\Delta_0,I)\cle 0$ and \r{csdilmi} is \r{lmisdi} for $P$ replaced by $-P$,
we first conclude that $I-D\Del_0$ is invertible. Hence, \r{fb} can be rewritten as \r{fbl} for $\nl(z)=\Del_0z$. As a consequence, this loop does have a response for all initial conditions $x_0\in\R^n$.
By the nominal stability assumption, all these responses satisfy $\lim_{t\to\infty}x(t)=0$. Hence,
$A+B\Del_0(I-D\Del_0)^{-1}C$ is Hurwitz and, thus, $X\cg 0$ by Theorem~\ref{Tlmipos}.

Due to \r{csdilmi} and again by Corollary~\ref{Csdi}, there exists some $\ve>0$ such that
\eql{h3}{
x(T)^T\Xs x(T)+\ve\int_{0}^T \|x(t)\|^2+\|u(t)\|^2\,dt+
\int_{0}^{T}
\mat{c}{y(t)\\u(t)}^T\!\!P\mat{c}{y(t)\\u(t)}
\,dt\leq x(0)^T\Xs x(0)
}
holds for all $T\geq 0$ and all admissible trajectories of \r{sy}.
Now take any $\nl$ satisfying \r{qccir} and any response of \r{fb}.
Then, \r{qccir} implies \eql{iqc0}{\int_0^T \mat{c}{z(t)\\w(t)}^T\!\!P\mat{c}{z(t)\\w(t)}\,dt\geq 0\te{for all}T\geq 0.}
Moreover, \r{h3} stays valid for this particular trajectory. Due to $y=z$ and $u=w$,
both inequalities are combined to infer \r{h0}. This implies \r{h0s} and thus \r{sta} for
$\ga:=\sqrt{\la_{\rm max}(\Xs)/\ve}$.

\subsection{Discussion of the circle criterion}

The exposition of dissipativity theory reveals that the circle criterion can be
formulated in various seemingly different but equivalent ways. For example, the assumption on $A$ in Theorem~\ref{Tcir} only plays a role in how to formulate the FDI. If using the general version of the FDI in Corollary~\ref{Csdi}, it remains valid for any linear system \r{sy} without the need for further adaptations. In view Corollary~\ref{Csdi} and Theorem~\ref{Tlmipos}, the hypotheses  in Theorem~\ref{Tcir} can be equivalently expressed through the existence of a solution of the dissipation inequality \r{csdilmi} that is positive definite. The criterion also admits an equivalent formulation in terms of an FDI in the closed right half-plane, which is not pursued in the current article.

The proof of Theorem~\ref{Tcir} only relies on the dissipation inequality \r{iqc0}
as long as $w(t)=\nl(z(t))$ holds for $t\geq 0$. Therefore, the result remains valid for multivalued static or dynamic mappings $\nl$ without the need for another proof (as exposed in detail in the section ``\nameref{Srsh}'').

\subsection[Combination of multipliers]{Combination of multipliers}
\label{Scom}

It was explicitly proposed in  \cite{MoyHil78,SafAth81} to work with whole families of supply rates to strengthen dissipativity-based stability tests. This leads to the following extension of the LMI formulation of the circle criterion, which requires no new proof.

\corollary{\label{Ccir}Let $\bm{P}\subset\S^{k+m}$ and $\nl$ be any nonlinearity satisfying \r{qccir} for all $P\in\bm{P}$.
If there exist some $P\in\bm{P}$ and $\Xs\cg0$ for which the dissipation LMI \r{csdilmi} holds,
then \r{fb} is stable.
}

The construction of such sets of supply rates for the Lur'e setup with a saturation is illustrated next. If the loop in Fig.~\ref{fig1} is
multivariable, the nonlinearity $\nl$ is  ``diagonally repeated'' as
\eql{sat}{
\nl(z_1,\ldots,z_\nw):=\col\mat{c}{\sat(z_1),\ldots,\sat(z_\nw)}.
}
Then, \r{qccir} is true for $P=P_{0_m,I_m}$, where $0_m$ and $I_m$ denote the zero and identity matrices of dimension $m$.
However, there are many more such choices.
If $e_j\in\R^\nw$ denotes the standard unit vector, then
$[e_j^Tz-e_j^T\nl(z)]^Te_j^T\nl(z)=(z_j-\sat(z_j))\sat(z_j)\geq 0$ for all $z\in\R^\nz$ and $j=1,\ldots,m$.
If multiplying these inequalities with $2\la_j\geq 0$ and summing, infer that \r{qccir} holds for \r{sat} and any
matrix $P$ in the convex conic hull
$$
\bm{P_1}:=
\left\{\sum_{j=1}^\nw \la_j\mat{cc}{0&e_je_j^T\\e_je_j^T&-2e_je_j^T}\st \la_1\geq 0,\ldots,\la_\nw\geq 0\right\}.
$$
Note that $\bm{P_1}$ is a convex set with an explicit LMI representation in terms of generators.
Convexity arguments lead to a larger  implicitly defined set $\bm{P_2}$, again with an LMI representation.

\lemma{Using $E(\delta_1,\ldots,\delta_\nw):=\diag(\delta_1,\ldots,\delta_\nw)$ let
$$
\bm{P_2}:=\left\{P=\mat{cc}{Q&S\\S^T&R}\in\S^{2\nw}\st R\cle 0,\
\mat{c}{I\\E(\delta)}^T\!\!
P
\mat{c}{I\\E(\delta)}
\cge 0\te{for all}\delta\in\{0,1\}^\nw\right\}.
$$
Then, \r{qccir} holds true for the nonlinearity \r{sat} and all $P\in\bm{P_2}$.
}
\proo{Select any $P\in\bm{P_2}$. Since $R\cle 0$, the mapping $\la\mapsto s(E(\la),I)$ is concave. Since
$[0,1]^\nw$ is the convex hull of $\{0,1\}^\nw$, infer from $P\in\bm{P_2}$ that $s(E(\la),I)\cge 0$ for all $\la\in[0,1]^\nw$.
Now choose $z\in\R^\nz$. Then, there exist $\la_j\in[0,1]$ with $\sat(z_j)=\la_jz_j$, for $j=1,\ldots,\nw$. This
implies $\nl(z)=E(\la)z$ and thus $s_P(\nl(z),z)=s_P(E(\la)z,z)=z^Ts_P(E(\la),I)z\geq 0$. 
}

These are just two examples of how to construct suitable classes of $\bm{P}$ that can be used in  Corollary~\ref{Ccir}.
The resulting tests for these classes are numerically verifiable. For $j=1$ or $j=2$, the containment $P\in\bm{P_j}$ can be expressed by a finite number of LMI constraints, such that $P\in\bm{P_j}$, $\Xs\cg 0$ and \r{csdilmi} form a system of LMIs whose feasibility can be numerically verified by any off-the-shelf solver. With a parser such as Yalmip \cite{Lof04}, it is a matter of minutes to implement such a test in Matlab \cite{ML20b}.

The numerical complexity of these tests depends on the complexity of the description of the sets $\bm{P_j}$ in terms of LMIs.
The cheapest is obtained with a singleton $\bm{P_0}:=\{P_{0_\nw,I_\nw}\}$. The complexity increases by taking $\bm{P_1}$, which involves $\nw$ extra linear programming inequalities, while this number grows to $2^\nw$ for $\bm{P_2}$. On the other hand, the test is least conservative for the latter choice. The article \cite{FetSch17a} provides a detailed discussion of how to generate further classes of supply rates for repeated sector-bounded nonlinearities in a systematic fashion.


\newpage
\section[Link to classical
integral
quadratic constraints theorem]
{Link to classical IQC theorem}
\label{sb:iqc}

We address the precise relation of Theorem~\ref{Trs} to the classical integral
quadratic constraints (IQC) theorem in \cite{MegRan97}.
For this purpose, we review some basic concepts in a concise fashion.

First, $\Le_2^n$ denotes the linear space of square integrable functions on $[0,\infty)$ equipped with the norm $\|.\|_2$
and $\hat x$ is the Fourier transform of $x\in\Le_2^n$.
The truncation $x_T$ of some signal $x:[0,\infty)\to\R^n$ for any $T\geq 0$ is defined by
$x_T=x$ on $[0,T]$ and $x_T=0$ on $(T,\infty)$. Then, $\Le_{2e}^n$ denotes the linear space of measurable
functions $x:[0,\infty)\to\R^n$ with $x_T\in\Le_2^n$ for all $T\geq 0$. Moreover, a map
$\Delta:\Le_2^\nz\to\Le_2^\nw$ is causal if $\Del(z)_T=\Del(z_T)_T$ holds for all
$T\geq 0$ and all $z\in\Le_{2e}^\nz$. Finally, $\Del$ is stable if there exists some $\ga\geq 0$ such that
$\|\Delta(z)_T\|_2\leq\ga\|z_T\|_2$ holds for all $z\in\Le_{2e}^\nz$ and all $T\geq 0$.

As depicted on the left in Fig.~\ref{Frs}, consider the feedback loop
\eql{fbi}{\dot x=Ax+Bw,\ \ z=Cx+Dw,\ \ w=\Del(z)+d,\te{and}x(0)=0,}
with an external disturbance input $d\in\Le_{2e}^\nw$. It is assumed that
$A$ is Hurwitz and $\Delta:\Le_{2e}^\nz\to\Le_{2e}^\nw$ is causal and stable.
This loop is well-posed if, for all $d\in\Le_{2e}^\nw$, there exists a unique response $z\in\Le_{2e}^\nw$ of \r{fbi} such that
the correspondingly defined map $d\mapsto z$ is causal. The loop is stable if, in addition, the map $d\mapsto z$ is stable.

For the rational multiplier defined by $\Pi(s):=\Psi(-s)^TP\Psi(s)$ with some $P\in\S^p$ and some stable filter $\Psi=\mat{cc}{\Psi_1&\Psi_2}$ of dimension $p\times(\nz+\nw)$, the IQC theorem in \cite{MegRan97} reads as follows. It is emphasized that this result can also be applied for irrational multipliers $\Pi$ and extends to a more abstract operator framework, as exposed in \cite{Meg01}, for example.

\theorem{\label{Tiqc}
Suppose that \r{fbi} is well-posed for $\tau \Del$ replacing $\Delta$ and any $\tau\in[0,1]$.
If $\Delta$ satisfies the frequency domain IQC
\eql{iqc}{
\int_{-\infty}^\infty
\mat{c}{\hat z(\io)\\ \tau\widehat{\Delta(z)}(\io) }^*\!\!\Psi(\io)^*P\Psi(\io)
\mat{c}{\hat z(\io)\\ \tau\widehat{\Delta(z)}(\io) }\,d\om\geq 0
\te{for all}z\in\Le_2^\nz,\ \tau\in[0,1]
}
and the FDI \r{fdis} is satisfied, then the loop \r{fbi} is stable.
}

With a minimal realization \r{psi} of $\Psi$, stability of $\Psi$ implies that the matrix $A_\Psi$ is Hurwitz. By Parseval's theorem, it is routine to verify that \r{iqc} for $\tau=1$ is
\eql{iqctd}{
\int_0^\infty v(t)^T\mi v(t)\,dt\geq 0
\te{for all} 
v=\Psi_1z+\Psi_2\Del(z)\te{with}z\in\Le_2^\nz.
}
On the other hand, the IQC \r{iqcZ} in Definition~\ref{DiqcZ} translates in the current setting into
\eql{iqcZ2}{
\int_0^T v(t)^T\mi v(t)\,dt-\xi(T)^TZ\xi(T)\geq 0\te{for all}T\geq0,\
v=\Psi_1z+\Psi_2\Delta(z)\text{\ with\ }z\in\Le_{2e}^\nz
}
with the signal space $\Le_{2e}^\nw\times \Le_{2e}^\nz$ replacing $\Se^\nw\times\Se^\nz$.
Following \cite{MegRan97}, \r{iqctd} is a soft dynamic IQC, while \r{iqcZ2} for $Z=0$ is a hard IQC. As emphasized in a footnote in \cite{MegRan97}, a soft IQC might be valid, while a hard IQC is not.
The following summary of results from \cite{SchVee18} clarifies that the concept of IQCs with a nontrivial terminal cost $Z\neq 0$ is adequate to encompass and generalize results
based on both hard and soft IQCs.

\remark{
For $\tau=0$, note that \r{iqc}  implies $\Psi_1^*P\Psi_1\cge 0$. Then,
it causes no loss of generality to assume that the multiplier satisfies
$\Psi_1^*P\Psi_1\cg 0$.
This is achieved by substituting $\Psi^*P\Psi$ with $\Psi_e^*P_\eps\Psi_e$ in Theorem \ref{Tiqc}, where
$$
\Psi_e:=\mat{cc}{\Psi_1&\Psi_2\\I_\nz&0},\
\mi_\ve:=\mat{cc}{\mi&0\\0&\ve I_\nz},\te{and $\ve>0$ is sufficiently small.}
$$
Since $\Psi_e^*P_\eps\Psi_e\cge \Psi^*P\Psi$, \r{iqc} persists to hold for $\Psi_e^*P_\eps\Psi_e$.
Moreover, the left-upper $k\times k$-block of $\Psi_e^*P_\eps\Psi_e$ equals $\Psi_1^*P\Psi_1+\eps I$ and hence satisfies $\Psi_1^*P\Psi_1+\eps I\cg 0$. Finally, since
$\col(\sy,I)^*\Psi_e^*\mi_\ve \Psi_e\col(\sy,I)=
\col(\sy,I)^*\Psi^*\mi \Psi\col(\sy,I)+\ve\sy^*\sy$ and $\sy$ is stable, \r{fdis} is valid for $\Psi_e^*P_\eps\Psi_e$ if $\ve>0$ is small enough.\epro}

\theorem{\label{Tiqc2}
Suppose that $\Psi_1^*P\Psi_1\cg 0$ and consider the following statements:
\enu{
\item The hypotheses in Theorem~\ref{Tiqc} are satisfied.
\item The hypotheses in Theorem~\ref{Tiqc} are satisfied for $\tau=1$ and $\Psi_2^*P\Psi_2\cle 0$ is valid.
\item The hypotheses in Theorem~\ref{Tiqc} are satisfied for $\tau=1$ and the loop \r{fbi} is stable.
\item
There exists some $Z=Z^T$ such that $\Del$ satisfies the IQC \r{iqcZ2} with terminal cost and
the dissipation LMI \r{lmif} has a solution $X=X^T$, which is coupled with $Z$ as in \r{pos}.
\item There exists a solution $Z=Z^T$ of the algebraic Riccati equation
\eql{are}{
A_\Psi^T Z+ ZA_\Psi-C_\Psi^T\mi C_\Psi+(ZB_\Psi-C_\Psi^T\mi D_\Psi)(D_\Psi^T\mi D_\Psi)^{-1}(B_\Psi^T Z-D_\Psi^T\mi C_\Psi)=0}
such that $A_\Psi+B_\Psi(D_\Psi^T\mi D_\Psi)^{-1}(B_\Psi^T Z-D_\Psi^T\mi C_\Psi)$ is Hurwitz.
}
Then a) or b) implies c). Moreover, c) implies d) and e), and $Z$ in d) can be taken from e).
}

This result is covered in \cite{SchVee18}, while related ones for multivalued uncertainties were developed in \cite{SchHol18}.
Observe that  a) $\Rightarrow$ c) restates Theorem~\ref{Tiqc}.
The implication b) $\Rightarrow$ c) goes back to \cite{Sei15,VeeSch13a}. It has the practically relevant consequence that there is no need to verify the assumptions in Theorem~\ref{Tiqc} for all homotopy parameters $\tau\in[0,1]$ to prove stability of \r{fbi}. Finally, the last statement is a combination of Lemma 4 and Lemma 7 in \cite{SchVee18} if recalling the change of sign in Definition~\ref{DiqcZ}.

As a consequence of Theorem~\ref{Tiqc2}, both a) or b) implies d) and, hence, the validity of the assumptions in Theorem~\ref{Trs}.
In this sense, Theorem~\ref{Trs} generalizes the  classical IQC theorem. This includes an explicit way to determine some
 suitable terminal cost matrix, as characterized in e).

If the loop \r{fbi} is assumed to be well-posed, then d) also implies c). Indeed,
stability of \r{fbi} is implied by d), as shown in the section ``\nameref{Siqc}.'' Moreover, \r{iqcZ2}
leads to \r{iqctd} by taking the limit $T\to\infty$, since $z\in\Le_{2}^\nz$ implies $\Delta(z)\in\Le_2^\nw$ by stability of $\Delta$ and, therefore, $v\in\Le_2^p$ and $\xi(T)\to 0$ for $T\to\infty$ (since $A_\Psi$ is Hurwitz).

Therefore, under all the assumptions in this section, c) and d) are essentially equivalent, which
substantiates the claim that Theorem~\ref{Trs} constitutes a seamless link and extension of robustness
results based on soft IQCs and dissipativity theory.

%

\newpage
\section[Zames-Falb multipliers]
{Zames-Falb multipliers}
\label{sb:zf}

Zames-Falb multipliers \cite{Osh67,WilBro68,ZamFal68} lead to much more powerful IQCs than the static ones, as used in the circle criterion. Note that convexity is key in proving a suitable dynamic IQC, even when the nonlinearities are multivalued.

{\bf Assumption.} $\vp=\partial f$ holds for some convex function $f:\R^\nz\to\R$ with $0\in\partial f(0)$.

Here, $\partial f$ denotes the multivalued map that takes $z\in\R^\nz$ into the set $\partial f(z)$ of all vectors $g\in\R^\nz$, for which the following subgradient inequality holds true:
\eql{subg}{
g^T(z-y)\geq f(z)-f(y)\te{for all}y\in\R^\nz.
}
Recall that $0\in\partial f(0)$ is equivalent to $f$ having a global minimum at $0\in\R^\nz$. As a consequence of the assumptions, $f$ is actually globally nonnegative. This leads to the following result.

\lemma{\label{Lzf}Let $h:[0,\infty)\to\R$ be taken with $h(t)\geq 0$ for $t\geq 0$ and $\int_0^\infty h(t)\,dt\leq 1$.
If $z\in\Se^\nz$ is filtered as $y(t):=z(t)-\int_0^t h(t-\tau)z(\tau)\,d\tau$ for $t\geq 0$ and if $w\in\Se^\nz$
is taken with $w(t)\in\partial f(z(t))$ for $t\geq 0$, then
\eql{subg2}{
\int_0^T w(t)^Ty(t)\,dt\geq 0\te{for all}T\geq 0.
}
}
\proo{Extend $z$ by $z(t):=0$ for $t<0$ to a mapping $z:\R\to\R^\nz$.
Fix $\tau\geq 0$. If using \r{subg} for $(g,z,y)$ replaced by $(w(t),z(t),z(t-\tau))$,
integration over  $t\in[0,T]$ gives
\eql{subg3}{
\int_0^T w(t)^T(z(t)-z(t-\tau))\,dt\geq \int_0^Tf(z(t))\,dt-\int_0^T f(z(t-\tau))\,dt\te{for all}T\geq 0.
}
The right-hand side is nonnegative: It equals $\int_0^Tf(z(t))\,dt-\int_{-\tau}^{T-\tau} f(z(t))\,dt$. For
$0\leq T<\tau$, this is $\int_0^Tf(z(t))\,dt\geq 0$; for $\tau\leq T$ we obtain
$\int_0^Tf(z(t))\,dt-\int_{0}^{T-\tau} f(z(t))\,dt=\int_{T-\tau}^Tf(z(t))\,dt\geq 0$.

Since $h(\tau)\geq 0$, we can multiply  \r{subg3} with $h(\tau)$ and integrate over $\tau\in[0,\infty)$ to infer
\eqn{
\int_0^\infty h(\tau)\,d\tau\int_0^T w(t)^Tz(t)\,dt-\int_0^Tw(t)^T\int_0^\infty h(\tau)z(t-\tau)\,d\tau\,dt\geq 0\te{for all}T\geq 0.
}
Since $1\geq \int_0^\infty h(\tau)\,d\tau\geq 0$, \r{subg3} for $\tau=0$ gives
$\int_0^T w(t)^Tz(t)\,dt\geq \int_0^\infty h(\tau)\,d\tau\int_0^T w(t)^Tz(t)\,dt$ for all $T\geq 0$.
Because of $\int_0^\infty h(\tau)z(t-\tau)\,d\tau=\int_0^t h(t-\tau)z(\tau)\,d\tau$, the proof is concluded.
}

Now suppose that $h$ in Lemma~\ref{Lzf} is the impulse response of some finite-dimensional
linear system with a minimal realization in terms of $(A_h,B_h,C_h,D_h)$. Then, $A_h$ is Hurwitz and $D_h=0$.
Moreover, \r{subg2} means $\int_0^T v(t)^TP_{\pas}v(t)\,dt\geq 0$ for all $T\geq 0$ along the trajectories of
\eql{zfm}{
\dot\xi=(I_k\otimes A_h)\xi+(I_k\otimes B_h)z,\ \xi(0)=0,\ v=\mat{c}{-I_k\otimes C_h\\0}\xi+\mat{c}{I_k\\0}z+\mat{c}{0\\I_k}w,\ w\in\partial f(z).
}
If $\Psi_{h}$ denotes the transfer matrix of the linear system in \r{zfm}, then this expresses
the fact that $\partial f$ satisfies an IQC with terminal cost matrix $0$ for the Zames-Falb multiplier
$\Psi_{h}^*P_{\pas}\Psi_{h}$.

\newpage
\section[Proof of Theorem~\ref{Tdy}]
{Proof of Theorem~\ref{Tdy}}
\label{sb:Tdy}

If considering the dissipation LMI \r{lmiK} for $\psi^*M\psi\cg0$, we introduce $\mat{cc}{\t C_\psi&\t D_\psi}$ with
\eql{h200}{
\mat{cc}{A_\psi&B_\psi\\I&0}^T\!\!\mat{ccc}{0&K\\K&0}\mat{cc}{A_\psi&B_\psi\\I&0}-
\mat{cc}{C_\psi&D_\psi}^T\!\!M\mat{cc}{C_\psi&D_\psi}+\mat{cc}{\t C_\psi&\t D_\psi}^T\!\!\mat{cc}{\t C_\psi&\t D_\psi}=0.
}
Then, $\mat{cc}{\t C_\psi&\t D_\psi}^T\mat{cc}{\t C_\psi&\t D_\psi}\cg 0$. With
$\t \psi(s):=\t C_\psi(sI-A_\psi)^{-1}B_\psi+\t D_\psi$, the FDI $\t\psi^*\t\psi\cg 0$ is trivially certified by $0$, and any sufficiently small $\t X\cg 0$ still certifies the same FDI.
If $\t \psi$ has $\t k$ rows, it was already shown that $\t\psi^*\t\psi=\t\psi^*(I_{\t k})\t\psi\cg 0$ implies
\eql{h201}{
0\cl \mat{c}{1\\\delta}^*\!\!P_0\mat{c}{1\\\delta}\t\psi^*\t\psi=
\mat{c}{\t\psi\\\delta\t\psi}^*\!\!(P_0\otimes I_{\t \nz})\mat{c}{\t\psi\\\delta\t\psi}.
}
Now note that \r{fdidel} has a certificate $X_\delta\cg 0$ by Corollary~\ref{Csdi}.
Due to Lemma~\ref{Lmul} and with the natural realization of $\col(\t\psi,\delta\t\psi)$,
\r{h201} hence admits the certificate $\diag(X_\delta,\ve\t X)$
for some small $\ve>0$, and this matrix is positive definite. By Lemma~\ref{Lmin}, \r{h201} still has a positive definite certificate if taking any minimal realization of $\col(\t\psi,\delta\t\psi)$. Now note that $\col(\t\psi,\delta\t\psi)=\col(\t\psi,\t\psi \delta)$, and recall that all minimal realizations of these two transfer matrices are related by a state-coordinate change.
Therefore, there also exists a positive definite certificate of
\eql{h202}{
0\cl
\mat{c}{\t\psi\\\t\psi \delta}^*\!\!(P_0\otimes I_{\t \nz})\mat{c}{\t\psi\\\t\psi \delta}=
\mat{c}{1\\\delta}^*\!\!\diag(\t\psi,\t\psi)^*(P_0\otimes I_{\t \nz})
\diag(\t\psi,\t\psi)\mat{c}{1\\\delta}
}
if taking a minimal realization of $\diag(\t\psi,\t\psi)\col(1,\delta)$. Instead, choose
the one obtained for the series interconnection of $\diag(\t\psi,\t\psi)$ and $\col(1,\delta)$ from
minimal realizations of the individual factors. Despite that this realization is possibly not minimal,
it is still stable. Again by Lemma~\ref{Lmin}, the corresponding dissipation LMI also admits a positive definite solution.
Hence, the related dissipation inequality leads to a hard IQC. Precisely, for $z\in\Se$ and $w=\delta z$,  the response of
\eql{h203}{
\dot\xi=\diag(A_\psi,A_\psi)\xi+\diag(B_\psi,B_\psi)\col(z,w),\ \xi(0)=0
}
and $\t v=\diag(\t C_\psi,\t C_\psi)\xi+\diag(\t D_\psi,\t D_\psi)\col(z,w)$ satisfies
\eql{h204}{
\int_0^T \t v(t)^T(P_0\otimes I_{\t \nz})\t v(t)\,dt\geq 0\te{for all}T\geq 0.
}
We now exploit that \r{h200} directly leads to
\mul{\arraycolsep.8ex
(\pl)^T\mat{ccc}{0&P_0\otimes K\\P_0\otimes K&0}
\mat{cc}{\diag(A_\psi,A_\psi)&\diag(B_\psi,B_\psi)\\I&0}-\\-\arraycolsep.8ex
(\pl)^T(P_0\otimes M)\mat{cc}{\diag(C_\psi,C_\psi)&\diag(D_\psi,D_\psi)}+\\+\arraycolsep.8ex
(\pl)^T(P_0\otimes I_{\t \nz})\mat{cc}{\diag(\t C_\psi,\t C_\psi)&\diag(\t D_\psi,\t D_\psi)}=0.
}
Therefore, by Lemma~\ref{Ldise}, the two outputs 
$v=\diag(C_\psi,C_\psi)\xi+\diag(D_\psi,D_\psi)\col(z,w)$ and
$\t v=\diag(\t C_\psi,\t C_\psi)\xi+\diag(\t D_\psi,\t D_\psi)\col(z,w)$
of \r{h203}
are related as
$$
\xi(T)^T(P_0\otimes K)\xi(T)-\int_0^T v(t)^T(P_0\otimes M)v(t)\,dt+\int_0^T \t v(t)^T(P_0\otimes I_{\t k})\t v(t)\,dt=0\te{for all}T\geq0.
$$
Then, \r{h204} guarantees that the claimed inequality holds true, which finishes the proof.
\epro

\newpage
\section[Auxiliary facts]
{Auxiliary facts}
\label{sb:aux}

\lemma{\label{Lschur} Let $P\in\S^{n_1+n_2}$ be a partitioned matrix with $P_1\cl 0$ and for $\al,\be>0$, consider
$$
\mat{cc}{P_{1}&P_{12}\\P_{12}^T&P_2}\cl 0,\ \mat{cc}{\al P_{1}&P_{12}\\P_{12}^T&P_2}\cl 0,\ \
\mat{cc}{\beta P_{1}&\beta P_{12}\\\beta P_{12}^T&P_2}\cl 0.
$$
The first inequality holds iff the Schur complement satisfies $P_2-P_{12}^TP_1^{-1}P_{12}\cl 0$. Let $P_2\cl 0$.
Then, there exists a large $\al>0$ (a small $\beta>0$) such that the second (third) inequality is valid. }

\proo{The first fact is standard \cite{BoyGha94}. The other statements follow, since the last two inequalities
are equivalent to $P_2-\frac{1}{\al}P_{12}^TP_1^{-1}P_{12}\cl 0$ and $P_2-\beta P_{12}^TP_1^{-1}P_{12}\cl 0$, respectively.
}

\lemma{\label{Lmul}Suppose that $\sy_1^*P\sy_1\cg 0$ and  $\sy_2^*\sy_2\cg 0$ and consider the realizations
$$
\sy_1=\mas{c|c}{A_1&B_1\hl C_1&D_1},\
\sy_2=\mas{c|c}{A_2&B_2\hl C_2&D_2},\
\sy_1\sy_2=\mas{cc|c}{A_1&B_1C_2&B_1D_2\\0&A_2&B_2\hl C_1&D_1C_2&D_1D_2}.
$$
If $\sy_1^*P\sy_1\cg 0$ and $\sy_2^*\sy_2\cg 0$ are certified by $X_1=X_1^T$ and $X_2=X_2^T$, then
$(\sy_1\sy_2)^*P(\sy_1\sy_2)\cg 0$ is certified by $\diag(X_1,\ve X_2)$ for some small $\ve>0$.
}

\proo{
The strict dissipation LMI involving $X_1$ can be perturbed with some $\ve>0$ to
$$
\mat{ccc}{A_1&B_1\\I&0}^T\!\!\mat{cc}{0&X_1\\X_1&0}\mat{ccc}{A_1&B_1\\I&0}-
\mat{ccc}{C_1&D_1\\0&I}^T\!\!\mat{cc}{P&0\\0&-\ve I}\mat{ccc}{C_1&D_1\\0&I}\cl0.
$$
Right-multiplying $\diag(I,(C_2\ D_2))$ and left-multiplying the transpose gives
\eql{d1}{\arraycolsep.4ex
(\pl)^T
\mat{cc|cc}{0&0&X_1&0\\0&0&0&0\hl X_1&0&0&0\\0&0&0&0}
\mat{cc|c}{A_1&B_1C_2&B_1D_2\\0&A_2&B_2\hl I&0&0\\0&I&0}-
(\pl)^T\mat{cc}{P&0\\0&-\ve I}
\mat{cc|c}{C_1&D_1C_2&D_1D_2\\0&C_2&D_2}
\cle 0.}
By trivially inflating the strict dissipation LMI for $X_2$ and multiplying with $\ve>0$, we infer
\eql{d2}{\arraycolsep.4ex
(\pl)^T
\mat{cc|cc}{0&0&0&0\\0&0&0&\ve X_{2}\hl 0&0&0&0\\0&\ve X_{2}&0&0}
\mat{cc|c}{A_1&B_1C_2&B_1D_2\\0&A_2&B_2\hl I&0&0\\0&I&0}-
(\pl)^T\mat{cc}{0&0\\0&\ve I}\mat{cc|c}{C_1&D_1C_2&D_1D_2\\0&C_2&D_2}\cle 0.}
If $x=\col(x_1,x_2,x_3)$ is a kernel vector of the left-hand side of \r{d2}, we infer $x_2=0$ and $x_3=0$; if $x$ is also a kernel vector of \r{d1}, then $x_1=0$. Summing \r{d1} and \r{d2} hence gives
$$
(\pl)^T
\mat{cc|cc}{0&0&X_1&0\\0&0&0&\ve X_2\hl X_1&0&0&0\\0&\ve X_2&0&0}
\mat{cc|c}{A_1&B_1C_2&B_1D_2\\0&A_2&B_2\hl I&0&0\\0&I&0}-
(\pl)^T\mat{cc}{P&0\\0&0}
\mat{cc|c}{C_1&D_1C_2&D_1D_2\\0&C_2&D_2}
\cl 0.
$$
A trivial simplification concludes the proof. }

\lemma{\label{Lmin}Let a linear system with transfer matrix $\sy$ be strictly $s_P$-dissipative.
Then, the dissipation LMI for a minimal realization of $\sy$ has a positive     definite solution iff this holds for any stabilizable and detectable realization of $\sy$.}

\proo{With any realization $\sy(s)=C(sI-A)^{-1}B+D$, the dissipation LMI reads as
\eql{disf}{\arraycolsep.4ex
\mat{cc|c}{A_1&A_{12}&B_1\\A_{21}&A_2&B_2\hl I&0&0\\0&I&0}^T\!\!
\mat{cc|cc}{0&0&X_{1}&X_{12}\\0&0&X_{12}^T&X_{2}\hl X_{11}&X_{12}&0&0\\X_{12}^T&X_{22}&0&0}
\mat{cc|c}{A_1&A_{12}&B_1\\A_{21}&A_2&B_2\hl I&0&0\\0&I&0}
-\mat{cc|c}{C_1&C_2&D}^T\!\!P\mat{cc|c}{C_1&C_2&D}\cl 0}
where the partition of the matrices is yet to be specified. Moreover, note that the existence of a positive definite solution of this dissipation LMI is invariant under any state-coordinate change applied to the state-space
realization.

{\bf Step 1.} Suppose $X\cg 0$ satisfies \r{disf}, and let $(A,B)$ be stabilizable but not controllable. Then, choose
state-coordinates, such that $A_1$ is Hurwitz, $A_{12}=0$, $B_1=0$, and $(A_2,B_2)$ is controllable. Canceling the first block-row and block-column in \r{disf} directly leads to
\eql{disc}{
\mat{ccc}{A_2&B_2\\ I&0}^T\!\!
\mat{cc}{0&X_{2}\\X_{2}&0}
\mat{ccc}{A_2&B_2\\ I&0}
-\mat{ccc}{C_2&D}^T\!\!P\mat{ccc}{C_2&D}\cl 0.}
Therefore, $X_2\cg 0$ certifies $s_P$-dissipativity for $\sy(s)=C_2(sI-A_2)^{-1}B_2+D$,
which is a controllable realization.

Conversely, let $X_2\cg 0$ satisfy \r{disc}. Then, take any $K\cg 0$ with $A_1^TK+KA_1\cl 0$
(which exists since $A_1$ is Hurwitz) and
choose $X_2$, $X_{12}:=0$ and $X_1:=\al K$ in \r{disf}. The left-hand side of \r{disf} has
$\al (A_1^TK+KA_1)$ as its left-upper block, all other blocks do not depend on $\al$, and the right-lower $2\times 2$ block is negative definite by \r{disc}. Therefore, \r{disf} is valid for all sufficiently large $\al>0$  by Lemma~\ref{Lschur}.
This proves the claim for moving between stabilizable and controllable realizations.

{\bf Step 2.}  Now suppose $X\cg 0$ satisfies \r{disf}, and let $(A,C)$ be detectable but not observable. Then, take state-coordinates such that $A_1$ is Hurwitz, $A_{21}=0$, $C_1=0$, and $(A_2,C_2)$ is observable.
Perform a further state-coordinate change with the transformation matrix
$$
\mat{ccc}{I&-X_1^{-1}X_{12}\\0&I},\te{and introduce the Schur complement}X_s:=X_2-X_{12}^TX_1^{-1}X_{12}\cg 0.$$
With $\bar A_{12}:=A_{12}+A_1X_1^{-1}X_{12}-X_1^{-1}X_{12}A_2$ and $\bar B_1:=B_1-X_1^{-1}X_{12}B_2$, this leads to
\eql{diso}{\arraycolsep.4ex
\mat{cc|c}{A_1&\bar A_{12}&\bar B_1\\0&A_2&B_2\hl I&0&0\\0&I&0}^T\!\!
\mat{cc|cc}{0&0&X_{1}&0\\0&0&0&X_{s}\hl X_{1}&0&0&0\\0&X_{s}&0&0}
\mat{cc|c}{A_1&\bar A_{12}&\bar B_1\\0&A_2&B_2\hl I&0&0\\0&I&0}
-\mat{cc|c}{0&C_2&D}^T\!\!P\mat{cc|c}{0&C_2&D}\cl 0.
}
Again by canceling the first row and column in \r{diso}, we infer
\r{disc} for $X_s$ replacing $X_2$.

Conversely, if $X_s\cg 0$ satisfies \r{disc}, select
$K\cg 0$ with $A_1^TK+KA_1\cl 0$ and choose $X_2:=X_s$, $X_{12}:=0$, and $X_1:=\beta K$ in \r{disf}.
With the blocks from \r{disc}, this reads as
\eql{dish}{
\mat{ccc}{
\be(A_1^TK+KA_1)&\be KA_{12}&\be KB_1\\
\be A_{12}^TK   &\r{disc}   &\r{disc}\\
\be B_1^TK      &\r{disc}   &\r{disc}}\cl 0.
}
By  \r{disc} and Lemma~\ref{Lschur}, \r{dish} is indeed satisfied for all sufficiently small $\be>0$.
This proves the claim for moving between detectable and observable realizations.

{\bf Step 3.} Finally, suppose that $\sy(s)=C(sI-A)^{-1}B+D$ is a stabilizable and detectable realization.
This can be reduced to some realization $\sy(s)=\hat C(sI-\hat A)^{-1}\hat B+D$ which is controllable, without loosing detectability. A second step reduces it to $\sy(s)=\tilde C(sI-\tilde A)^{-1}\tilde B+D$, a controllable and observable realization. Then the dissipation LMI for $(A,B,C,D)$ has a positive definite solution iff this holds
for $(\hat A,\hat B,\hat C,D)$ (by Step 1.) iff this is true for $(\tilde A,\tilde B,\tilde C,D)$ (by Step 2.). This concludes  the proof.
}

\lemma{\label{Ldise}
For the symmetric matrices $P$, $\t P$,
let $K=K^T$ satisfy
$$
\mat{cc}{A_\Psi&B_\Psi\\I&0}^T\!\!\mat{ccc}{0&K\\K&0}\mat{cc}{A_\Psi&B_\Psi\\I&0}-
\mat{cc}{C_\Psi&D_\Psi}^T\!\!P\mat{cc}{C_\Psi&D_\Psi}+\mat{cc}{\t C_\Psi&\t D_\Psi}^T\t P\mat{cc}{\t C_\Psi&\t D_\Psi}=0.
$$
Then, the lossless dissipation ``inequality''
$$
\xi(T)^TK\xi(T)-\int_0^T v(t)^TPv(t)\,dt+\int_0^T \t v(t)^T\t P\t v(t)\,dt=\xi(0)^TK\xi(0)\te{for all}T\geq0
$$
holds along all trajectories of $\dot \xi =A_\Psi \xi+B_\Psi u$, $v=C_\Psi x+D_\Psi u$, $\t v=\t C_\Psi x+\t D_\Psi u$.
}
\proo{If right-multiplying the matrix equation with $\col(\xi(t),u(t))$ and left-multiplying the transpose, infer
$\frac{d}{dt} \xi(t)^TK\xi(t)-v(t)^TPv(t)+\t v(t)^T\t P\t v(t)=0$ for all $t\geq 0$.
Integration over $[0,T]$ for $T\geq 0$ completes the proof.
}

\processdelayedfloats 
\clearpage

\section{Author Biography}

Carsten W. Scherer received his Ph.D. degree in mathematics from the University of W\"uq
rzburg (Germany) in 1991. In 1993, he joined Delft University of Technology (The Netherlands), where he held positions as an assistant and associate professor. From 2001 to 2010, he was a full professor in the Delft Center for Systems and Control at Delft University of Technology. Since 2010, he has held the SimTech Chair for Mathematical Systems Theory in the Department of Mathematics at the University of Stuttgart (Germany).

Dr. Scherer served as the chair of the IFAC Technical Committee on Robust Control (2002-2008), and as associate editor for the IEEE Transactions on Automatic Control, Automatica, Systems and Control Letters, and the European Journal of Control. Since 2013 he is an IEEE Fellow ``for contributions to optimization-based robust controller synthesis.''

Dr. Scherer's main research activities cover various directions in applying optimization techniques for developing new advanced robust controller design algorithms and their application to mechatronics and aerospace systems.

{\bf Acknowledgements.} This work was funded by Deutsche Forschungsgemeinschaft (DFG, German Research Foundation) under Germany's Excellence Strategy - EXC 2075 - 390740016. We acknowledge the support by the Stuttgart Center for Simulation Science (SimTech).
Moreover, the author would like to sincerely thank Joost Veenman, Matthias Fetzer, and
Tobias Holicki for the joint work on the topics of this article. Many thanks as well to Siep Weiland for the long-term collaboration on our lecture notes about linear matrix inequalities in control.


\end{document}